\documentclass[titlepage,11pt]{article}
\usepackage{latexsym}
\usepackage{amsfonts}
\usepackage{amsmath}
\usepackage{tikz}
\usetikzlibrary{decorations.pathreplacing,calligraphy}
\newcommand\blackslug{\hbox{\hskip 1pt \vrule width 4pt height 8pt depth 1.5pt
        \hskip 1pt}}
\newcommand\bbox{\hfill \quad \blackslug \bigbreak}

\def\LL{,\ldots,}
\def\cupcup{\cup\cdots\cup}
\newcommand{\vare}{\varepsilon}

\title{Pure pairs. VIII. Excluding a sparse graph}
\author{Alex Scott\thanks{Research supported by EPSRC grant EP/V007327/1.}\\
Mathematical Institute, University of Oxford, Oxford OX2 6GG, UK
\\
\\
Paul Seymour\thanks{Supported by AFOSR grants A9550-19-1-0187 and FA9550-22-1-0234, and NSF
grants DMS-1800053 and DMS-2154169.}\\
Princeton University, Princeton, NJ 08544, USA
\\
\\
Sophie Spirkl\thanks{We acknowledge the support of the Natural Sciences and Engineering Research
Council of Canada (NSERC) [funding reference number RGPIN-2020-03912]. Cette
recherche a \'et\'e financ\'ee par le Conseil de recherches en sciences naturelles et
en génie du Canada (CRSNG) [num\'ero de r\'ef\'erence RGPIN-2020-03912].}\\
University of Waterloo, Waterloo, Ontario N2L3G1, Canada}

\date{July 21, 2020; revised \today}

\newtheorem{thm}{}[section]

\newcommand{\Proof}{\noindent{\bf Proof.}\ \ }

\begin{document}
\maketitle
\begin{abstract}
A pure pair of size $t$ in a graph $G$ is a pair $A,B$ of disjoint subsets of $V(G)$, each of cardinality at least $t$, such that $A$ is either complete or anticomplete to $B$.
It is known that, for every forest $H$, every graph on $n\ge2$ vertices that does not contain $H$ or its complement as an induced subgraph has 
a pure pair of size $\Omega(n)$; furthermore, this only holds when $H$ or its complement is a forest.

In this paper, we look at pure pairs of size $n^{1-c}$, where $0<c<1$.
Let $H$ be a graph: does every graph on $n\ge2$ vertices that does not contain 
$H$ or its complement as an induced subgraph have a pure pair of size $\Omega(|G|^{1-c})$?
The answer is related to the {\em congestion} of $H$,
the maximum of $1-(|J|-1)/|E(J)|$ over all subgraphs $J$ of $H$  with an edge. (Congestion is nonnegative,
and equals zero exactly when $H$ is a forest.)
Let $d$ be the smaller of the congestions of $H$ and $\overline{H}$. We show that the answer to the question above is ``yes'' 
if $d\le c/(9+15c)$,
and ``no'' if $d>c$.

\end{abstract}

\section{Introduction}
Graphs in this paper are finite, and without loops or parallel edges. Let $A,B\subseteq V(G)$ be disjoint.
We say that $A$ is {\em complete} to $B$, or $A,B$ are {\em complete}, if every vertex in $A$ is adjacent to every vertex in $B$,
and similarly $A,B$ are {\em anticomplete} if no vertex in $A$ has a neighbour in $B$. We say $A$ {\em covers} $B$ if every vertex
in $B$ has a neighbour in $A$. A {\em pure pair} in $G$ is a pair $A,B$ of disjoint subsets of $V(G)$ such that $A,B$
are complete or anticomplete. The number of vertices of $G$ is denoted by $|G|$.
The complement graph of $G$ is denoted by $\overline{G}$.
Let us say $G$ {\em contains} $H$ if some induced subgraph of $G$ is isomorphic to $H$,
and $G$ is  {\em $H$-free} otherwise.
If $X\subseteq V(G)$, $G[X]$ denotes the subgraph induced on $X$.

When can we guarantee that a graph has a large pure pair?  The most we can ask for is a pure pair $A,B$, where both sets have size 
$\Omega(|G|)$.  It turns out that to get this it is enough to exclude a forest and its complement.  In an earlier paper with Maria Chudnovsky, we proved the following~\cite{pure1}:
\begin{thm}\label{pure1}
For every forest $H$ there exists $\vare>0$ such that for every graph $G$ with $|G|>1$ that is both $H$-free and $\overline{H}$-free,
there is a pure pair $A,B$ in $G$ with $|A|,|B|\ge \vare|G|$.
\end{thm}
It is easy to see (with a random construction) that this has a converse: if $H$ is a graph and neither $H$ nor $\overline{H}$ is a forest then there is
no $\vare$ as in \ref{pure1}. 

What happens if $H$ is not a forest?  Do we still get a large pure pair?
A theorem of Erd\H{o}s, Hajnal and Pach~\cite{EHP} (not normally stated in this form, but this is equivalent) says that we do:
\begin{thm}\label{EHP}
For every graph $H$ there exist $\vare,b>0$ such that for every graph $G$ with $|G|>1$ that is both $H$-free and $\overline{H}$-free,
there is a pure pair $A,B$ in $G$ with $|A|,|B|\ge \vare|G|^b$.
\end{thm}
But in this, $b$ might be very small, and that raises the question: what sort of graph $H$ will make \ref{EHP} true with
$b$ close to $1$? 

In an earlier paper~\cite{pure5} we proved that certain sparse graphs $H$ have this property.
We say that $H$ has {\em branch-length} at least $k$ if $H$ is an induced subgraph of a graph that can be constructed as follows:
start with a multigraph (possible with loops or parallel edges) and subdivide each edge at least $k-1$ times.
We proved in~\cite{pure5} that:
\begin{thm}\label{stablehandle}
Let $c>0$ with $1/c$ an integer, and let $H$ be a graph with branch-length at least $4c^{-1}+5$. Then there exists $\vare>0$ such that
for every graph $G$ with $|G|>1$ that is
both $H$-free and $\overline{H}$-free, there is
a pure pair $A,B$ in $G$ with $|A|\ge \vare|G|$ and $|B|\ge \vare|G|^{1-c}$.
\end{thm}
So graphs $H$ with large branch-length make \ref{EHP} true with
$b$ close to $1$, but there are graphs with small branch-length that do this as well, for instance forests; so large branch-length
is sufficient but not necessary for our property. 

We also proved (unpublished) a similar theorem, that if $H$ can be obtained from a multigraph $H'$ by selecting a spanning tree
and subdividing many times all edges not in the tree,
then something like \ref{stablehandle} holds. (Not quite the analogue of \ref{stablehandle}: for a given value of $c$, the number of times the non-tree
 edges have to be subdivided depends not only on $c$ but also on the graph $H'$.) But that is not the answer either; we shall see, 
for instance, that 
if we take a long cycle, and for each of its
vertices $v$ add a new vertex adjacent only to $v$, this graph $H$ has our property, but cannot be built by either of
the constructions just given.

The answer is related to ``congestion''.
Let $H$ be a graph.
If $E(H)\ne \emptyset$, we define the {\em congestion}
of $H$ to be the maximum of $1-(|J|-1)/|E(J)|$, taken over all subgraphs
$J$ of $H$ with at least one edge; and if $E(H)=\emptyset$, we define
the congestion of $H$ to be zero. Thus the congestion of $H$ is always non-negative, and equals zero
if and only if $H$ is a forest. Graphs of small congestion must have large girth, but that is not the same thing: 
for instance, there are graphs with girth and average degree at least 100, and their congestion  is at least .98.
Roughly speaking, a graph has small congestion and only if it has large girth and its maximum average degree is at most slightly more than two
(that is, every induced subgraph has average degree at most $2+\vare$ for some small $\vare$). As we shall see, another way to think 
of graphs with 
small congestion is, they are the graphs that can be built by starting from the null graph and repeated adding vertices with at most 
one neighbour, and adding long paths joining vertices in what we already have built.

It turns out that graphs $H$ where one of $H,\overline{H}$ has small congestion satisfy \ref{EHP} with a value of $b$
close to $1$, while those where both of $H,\overline{H}$ have large congestion do not. Let us say these two things more precisely.
The first of these statements is the main result of the paper, the following:
\begin{thm}\label{mainthm}
Let $c>0$, and let $H$ be a graph such that one of $H,\overline{H}$ has congestion at most $\frac{c}{9+15c}$. 
Then there exists $\vare>0$ such that
for every graph $G$ with $|G|>1$ that is both $H$-free and $\overline{H}$-free, there is
a pure pair $A,B$ in $G$ with $|A|,|B|\ge \vare|G|^{1-c}$.
\end{thm}

The second statement is the following, which we will prove now:
\begin{thm}\label{conjconverse}
Let $c>0$, and let $H$ be a graph such that $H, \overline{H}$ both have congestion more than $c$.
There is no $\vare>0$ such that
for every graph $G$ with $|G|>1$ that is both $H$-free and $\overline{H}$-free, there is
a pure pair $A,B$ in $G$ with $|A|, |B|\ge \vare|G|^{1-c}$.
\end{thm}
\Proof Let $J$ be a subgraph of $H$ with $E(J)\ne \emptyset$ and $|J|-1< (1-c)|E(J)|$, and let $J'$ be a 
subgraph of $\overline{H}$ with $E(J')\ne \emptyset$ and 
$|J'|-1< (1-c)|E(J')|$. Let 
$$c':= 1-\max\left(\frac{|J|-1}{|E(J)|},\frac{|J'|-1}{/|E(J')|}\right);$$ 
so $c<c'<1$. Choose $d$ with $c<d<c'$.
Let $\vare>0$, let  $n$ be a large number, let $p:= n^{d-1}$, and let $G$ be a random graph on $n$ vertices, 
in which every pair of 
vertices are adjacent independently with probability $p$. Then (if $n$ is sufficiently large with $\vare$ given), an easy calculation
(which we omit) shows that, $G$ has no pure pair $A,B$ in $G$ with 
$|A|, |B|\ge (\vare/2) n^{1-c}$ with probability more than $1/2$ (indeed, approaching $1$ as $n$ goes to infinity). 

The expected number of induced subgraphs 
of $G$ isomorphic to $J$ is at most 
$$n^{|J|}p^{|E(J)|}= n^{|J|+(d-1)|E(J)|}\le n^{|J|-(|J|-1)\frac{1-d}{1-c'}}= n^{1-(|J|-1)\frac{c'-d}{1-c'}}\le  n/16$$
since $|E(J)|\ge (|J|-1)/(1-c')$. Consequently the probability that there are more than $n/4$ such subgraphs is at most $1/4$.
Similarly the probability that there are more than $n/4$ induced subgraphs isomorphic to $J'$ is at most $1/4$; and so
with positive probability, $G$ contains at most $n/4$ of copies of $J$, and most $n/4$ copies of $J'$, and has no 
pure pair $A,B$ in $G$ with
$|A|, |B|\ge (\vare/2)|G|^{1-c}$. But then by deleting at most $n/2$ vertices, we obtain a graph $G'$ containing neither $J$ 
nor $J'$, and hence containing neither $H$ nor $\overline{H}$, and with no pure pair $A,B$ with $|A|, |B|\ge \vare |G'|^{1-c}$ (since $\vare |G'|^{1-c}\ge  (\vare/2) n^{1-c}$).
This proves \ref{conjconverse}.~\bbox

The conclusion of \ref{stablehandle} is stronger than that of \ref{mainthm}: one of the sets of the pure pair has linear size.
That raises the question, is the corresponding strengthening of \ref{mainthm} true? More exactly:
\begin{thm}\label{conj}
{\bf Possibility:} For all $c>0$, there exists $\xi>0$ with the following property. For every graph $H$ with congestion
at most $\xi$, there exists $\vare>0$ such that for every graph $G$ with $|G|>1$
that is $H$-free and $\overline{H}$-free, there is
a pure pair $A,B$ in $G$ with $|A|\ge \vare|G|$ and $|B|\ge \vare|G|^{1-c}$.
\end{thm}
(The difference from \ref{mainthm} is that we are now asking for $|A|$ to be linear.) We were unable to decide this.

\section{Reduction to the sparse case}

Let us say a graph $G$ is {\em $\vare$-sparse} if every vertex has degree less than $\vare|G|$. An {\em anticomplete pair}
in $G$ is a pair $A,B$ of subsets of $V(G)$ that are anticomplete. For $\gamma,\delta\ge 0$, let us say $G$ is 
{\em $(\gamma,\delta)$-coherent} if there is no anticomplete pair $A,B$ with $|A|\ge \gamma$ and $|B|\ge \delta$.
We observe:
\begin{thm}\label{big}
If $\vare>0$, and $\vare\le 1/2$, and $G$ is $\vare$-sparse and $(\vare|G|,\vare|G|)$-coherent with $|G|>1$,
then $|G|> 1/\vare$.
\end{thm}
\Proof
Suppose that $|G|\le 1/\vare$.  If some distinct $u,v\in V(G)$ are non-adjacent, 
$\{u\},\{v\}$ form an anticomplete pair, both of cardinality at least $\vare|G|$, a contradiction.
So $G$ is a complete graph; but its maximum degree is less than $\vare|G|$ and $\vare\le 1/2$, which is impossible since $|G|>1$.
This proves \ref{big}.~\bbox

\bigskip

If $G$ is a graph and $v\in V(G)$, a {\em $G$-neighbour} of $v$ means a vertex of $G$ adjacent to $v$ in $G$.
A theorem of R\"odl~\cite{rodl} implies the following:
\begin{thm}\label{rodlthm}
For every graph $H$ and all $\eta>0$ there exists $\delta>0$ with the following property.
Let $G$ be an $H$-free graph. Then there exists $X\subseteq V(G)$ with $|X|\ge \delta |G|$, such that 
one of $G[X]$, $\overline{G}[X]$ is $\eta$-sparse.
\end{thm}

Consequently, in order to prove \ref{mainthm}, it suffices to prove the following:

\begin{thm}\label{sparsethm}
Let $c>0$, and let $H$ be a graph with congestion at most $\frac{c}{9 +15c}$. Then there 
exists $\vare>0$ such that 
every $\vare$-sparse $(\vare|G|^{1-c}, \vare|G|^{1-c})$-coherent graph $G$ with $|G|>1$ contains $H$.
\end{thm}
\noindent{\bf Proof of \ref{mainthm}, assuming \ref{sparsethm}.\ \ }
Let $c>0$, and let $H$ have congestion at most $\frac{c}{9 +15c}$.
Choose $\eta\le 1/2$ such that \ref{sparsethm} holds with $\vare$ replaced by $\eta$.
Choose $\delta$ such that \ref{rodlthm} holds. Let $\vare:=\eta\delta$. We claim that $\vare$ satisfies \ref{mainthm}.

Let $G$ be a graph with $|G|>1$ that is $H$-free and $\overline{H}$-free.
We must show that there is
a pure pair $A,B$ in $G$ with $|A|,|B|\ge \vare|G|^{1-c}$.
From the choice of $\delta$, there exists $X\subseteq V(G)$ with $|X|\ge \delta |G|$, such that        
one of $G[X]$, $\overline{G}[X]$ is $\eta$-sparse; and by \ref{big} we may assume that $|G|> 1/\vare\ge 1/\delta$, and so $|X|>1$.
If $G[X]$ is $\eta$-sparse, then from the choice of $\eta$, \ref{sparsethm} applied to $G[X]$ 
implies that there is
an anticomplete pair $A,B$ in $G[X]$ with 
$$|A|,|B|\ge \eta|X|^{1-c}\ge \eta\delta^{1-c}|G|^{1-c}\ge \eta\delta|G|^{1-c}=\vare|G|^{1-c},$$
as required. If $\overline{G}[X]$ is $\eta$-sparse we argue similarly, working in $\overline{G}[X]$. This proves \ref{mainthm}.~\bbox

The remainder of the paper is devoted to proving \ref{sparsethm}.

\section{Congestion}

In this section we replace the ``small congestion'' hypothesis of \ref{sparsethm} with a different hypothesis that is easier to use.
A {\em branch} of a graph $G$ is either:
\begin{itemize}
\item a path $P$ of $G$ with ends $p_1,p_2$ say, with length at least one, 
such that all the internal vertices of $P$ have degree two in $G$,
and $p_1,p_2$ have degree different from two in $G$; or
\item a cycle of $G$ such that all its vertices except at most one have degree two in $G$.
\end{itemize}
It follows that every edge of $G$ belongs to a unique branch of $G$.

We need two ways to make a larger graph from a smaller one. First, let $H$ be a graph, and let $v\in V(H)$ have degree 
at most one; then we say that $H$ is obtained from $H\setminus \{v\}$ by {\em adding the subleaf $v$}.
Second, let $H$ be a graph, and let $P$ be an induced path of $H$ of length at least two, such that all its internal vertices 
have degree two in $H$; then we say that $H$ is obtained from $H\setminus P^*$ by {\em adding the handle $P$}, where $P^*$
denotes the set of internal vertices of $P$. 

Let $\beta\ge 2$ be an integer. We say that a graph $H$ is {\em weakly $\beta$-buildable} if it can be constructed, 
starting from the null graph, by repeatedly
either adding a subleaf, or adding a handle of length at least $\beta$.
It is easy to see that if $H$ is weakly $\beta$-buildable,
then $H$ has congestion at most $1/\beta$.
We need a partial converse to this:
\begin{thm}\label{longbranch}
Let $\xi\le 1/3$, and $\beta=\lfloor 1/(3\xi)\rfloor +1$. If $H$ is non-null and has congestion at most $\xi$, then $H$ is weakly $\beta$-buildable.
\end{thm}
\Proof
We proceed by induction on $|H|$, and so we may assume that every vertex has degree at least two, and $H$
is connected.
If $C$ is an induced cycle of $H$, then since $H$, and therefore $C$, has congestion at most $\xi$, it follows that
$\xi\ge 1-(|C|-1)/|E(C)|$, and since $|C|=|E(C)|$, we deduce that $\xi\ge 1/|C|$, that is, $|C|\ge 1/\xi$. Thus 
every induced cycle, and hence every cycle, of $H$
has length at least 
$1/\xi>  1/(3\xi)+1$.
Suppose that some branch $B$ is a cycle.
Since $B$ has length more than  $1/(3\xi)+1$, $H$ can be obtained from a smaller graph by adding a handle of 
length more than $1/(3\xi)$ and hence at least $\beta$, and the result follows from the inductive hypothesis.

So we may assume that
every branch is a path with distinct ends, both in $W$, where
$W$ is the set of vertices of $H$ with degree at least three in $H$.
Thus every branch, of length $b$ say,
contains exactly $b-1$ vertices of degree two, and they belong to no other branches.

Let $H$ have $k$ branches, 
with lengths $b_1\LL b_k$ respectively.
Thus $H$ has $b_1+\cdots +b_k$ edges. 
At most two edges in each branch are incident with vertices in $W$;
and by summing the degrees 
of the vertices in $W$, we deduce that $|W|\le 2k/3$. Hence  
$$|H|=|W|+(b_1-1)+\cdots+(b_k-1)\le 2k/3+(b_1+\cdots+b_k)-k.$$
Since $1-(|H|-1)/|E(H)|\le \xi$, it follows that
$|H|\ge 1+(1-\xi)|E(H)|$, and so
$$2k/3+(b_1+\cdots+b_k)-k\ge 1+(1-\xi)(b_1+\cdots+b_k),$$ that is,
$b_1+\cdots+b_k\ge 1/\xi+k/(3\xi)$. Consequently some $b_i> 1/(3\xi)$, and in particular, the branch of maximum length
has length more than $1/(3\xi)$. Since $\xi\le 1/3$, this branch has length at least two, and so
$H$ can be obtained from a smaller graph by adding a handle of
length more than $1/(3\xi)$, and hence at least $\beta$, as required.
This proves \ref{longbranch}.~\bbox

Let $\beta\ge 2$ be an integer. We say that 
$G$ is {\em $\beta$-buildable} if it can be constructed, starting from a two-vertex graph with no edges, by repeatedly
adding a handle of length at least $\beta$. 
We observe:
\begin{thm}\label{weakbuild}
For $\beta\ge 2$, if $H$ is weakly $\beta$-buildable, then $H$ is an induced subgraph of a $\beta$-buildable graph.
\end{thm}
\Proof
We proceed by induction on $|H|$. We may assume that $H$ can be obtained from a weakly $\beta$-buildable graph $H'$
by either adding a subleaf,
or adding a handle of length at least $\beta$.
From the inductive hypothesis, $H'$ is an induced subgraph of a $\beta$-buildable
graph $J'$; and so $H$ is an induced subgraph of a graph $J$, where $J$ is obtained from $J'$ by either
adding a subleaf,
or adding a handle of length at least $\beta$. In the second case, $J$ is $\beta$-buildable, so we assume that $J$
is obtained from $J'$ by adding a subleaf $v$.
Thus $v$ has at most one neighbour in $V(J')$; and since $|J'|\ge 2$, 
we can add a handle
$B$ to $J'$ of length at least $\max(4,\beta)$, such that $v$ is an internal vertex of $B$. Consequently
$J$, and hence $H$, is an induced subgraph of a $\beta$-buildable graph
as required.  This proves \ref{weakbuild}.~\bbox

In order to prove \ref{sparsethm}, it therefore suffices to show the following (because \ref{sparsethm} is trivially 
true when $c\ge 1$, and if a graph has congestion at most $\frac{c}{9 +15c}$ when $c<1$ 
then by \ref{longbranch} and \ref{weakbuild} it is 
$\beta$-buildable with $c>1/\lfloor(\beta-3)/3\rfloor$):
\begin{thm}\label{handlethm}
Let $\beta\ge 2$ be an integer, let $H$ be a $\beta$-buildable graph, and let $c>1/\lfloor(\beta-3)/3\rfloor$.
There
exists $\vare>0$ such that
every $\vare$-sparse $(\vare|G|^{1-c}, \vare|G|^{1-c})$-coherent graph $G$ with $|G|>1$ contains $H$.
\end{thm}
This will be proved in the final section.

\section{Blockades, and a proof sketch}

Let $G$ be a graph and let the sets $B_i\;(i\in I)$ be nonempty, pairwise disjoint subsets of $V(G)$, where $I$
is a set of integers. We call $(B_i:i\in I)$ a {\em blockade} in $G$, and the sets $B_i\;(i\in I)$ are its {\em blocks}; 
its {\em length} is $|I|$, and its {\em width} is
$\min(|B_i|:i\in I)$. The {\em shrinkage} of $\mathcal{B}$ is the number $\sigma$ such that the width is $|G|^{1-\sigma}$.
(We will not need shrinkage until the end of the paper.)
If $\mathcal{B}=(B_i:i\in I)$ is a blockade in $G$, an induced subgraph $J$ of $G$
is {\em $\mathcal{B}$-rainbow} if $V(J)\subseteq \bigcup_{i\in I}B_i$, and each block of $\mathcal{B}$ contains
at most one vertex of $J$. If $H$ is a graph, a {\em $\mathcal{B}$-rainbow copy of $H$} means a $\mathcal{B}$-rainbow induced subgraph 
of $G$ that is isomorphic to $H$.

If $A,B$ are disjoint subsets of a graph $G$, the {\em max-degree} from $A$ to $B$ is defined to be the maximum, over $v\in A$,
of the number of neighbours of $v$ in $B$.
Let $\mathcal{B}=(B_i:i\in I)$ be a blockade in a graph $G$. For all distinct $i,j\in I$, let $d_{i,j}$ be the max-degree
from $B_i$ to $B_j$. The {\em linkage} of $\mathcal{B}$ is the maximum of $d_{i,j}/|B_j|$, over all distinct
$i,j\in I$ (or zero, if $|I|\le 1$).

We will prove in \ref{buildableforce} that if $H$ is a $\beta$-buildable graph, and $c>1/\lfloor(\beta-3)/3\rfloor$, and $G$ is $(|G|^{1-c}, |G|^{1-c})$-coherent and sufficiently large, then
there is an $\mathcal{A}$-rainbow copy of $H$ for every blockade $\mathcal{A}$ in $G$ with sufficient length and sufficiently small
shrinkage and linkage. This will imply \ref{handlethm}.

As we said earlier, we do not know whether
\ref{handlethm} is true with ``$ (\vare |G|^{1-c},\vare |G|^{1-c})$-coherent'' replaced by ``$ (\vare |G|^{1-c},\vare|G|)$-coherent''.
But the latter is sufficient for almost all the proof, and so we have written the proof just using this 
where we can.

The idea of the proof of \ref{buildableforce} is to work by induction on $|H|$; so we can assume that $H$ is obtained by adding 
a handle of length at least $\beta$ to a graph $H'$ for which the theorem holds. 
Since the theorem holds for $H'$, there are numbers 
$K',\lambda',\sigma'$, such that in every sufficiently large $(|G|^{1-c},|G|^{1-c})$-coherent graph $G$, and for every 
blockade $\mathcal{A}$ in $G$ with length at least $K'$ and with linkage and shrinkage at most $\lambda', \sigma'$ respectively, 
there must be an $\mathcal{A}$-rainbow copy of $H'$. In fact we prove more, that for all sufficiently small $\sigma'$ (less than
$c-1/\lfloor(\beta-3)/3\rfloor$) 
there exist $K', \lambda'$ with this property. We would like to prove the same statement for $H$, modifying $K,\sigma, \lambda$
appropriately. So we are given $\sigma$, less than $c-1/\lfloor(\beta-3)/3\rfloor$ (and it is important that this is a 
strict inequality).
We choose $\sigma'$ strictly between
$\sigma$ and $c-1/\lfloor(\beta-3)/3\rfloor$; now apply the inductive hypothesis to $H'$ and $\sigma'$ to get $K', \lambda'$; and use these
to construct $K,\lambda$ with the property we want. The rest of the paper is explaining the details of this last sentence, but
the crucial thing is that we are proving that $\sigma$ works for $H$, with the knowledge that a strictly larger number, $\sigma'$,
works for $H'$. This allows some wiggle room, which would not be available if we were trying to prove \ref{conj}.

Let us give some idea of the details just mentioned. The number $c$ is fixed throughout, and is larger than $1/\lfloor(\beta-3)/3\rfloor$. Let us say that
$(N,K,\sigma,\lambda,c)$ {\em forces} $H$ if for every $(|G|^{1-c},|G|^{1-c})$-coherent graph $G$ with $|G|\ge N$,
and for every blockade $\mathcal{A}$ in $G$ of length $K$,
shrinkage at most $\sigma$, and linkage at most $\lambda$, there is an $\mathcal{A}$-rainbow copy of $H$. 

We know that for all $\sigma'< c-1/\lfloor(\beta-3)/3\rfloor$, there exist $N',K',\lambda'$, such that 
$(N',K',\sigma',\lambda',c)$ forces $H'$; and we need to show that for all $\sigma< c-1/\lfloor(\beta-3)/3\rfloor$, there exist $N,K,\lambda$, such that 
$(N,K,\sigma,\lambda,c)$ forces $H$.
To obtain $H$ from $H'$, we need to add a handle of 
some specified length, at least $\beta$, with specified ends $u',v'$, but we can do this in two steps: first, add two vertices $u,v$ adjacent 
only to $u', v'$ respectively (forming $H''$ say); and then add a handle to $H''$ with ends $u,v$ of the right length 
(at least $\beta-2$).

The first part, adding the leaves $u,v$, is easy, by means of a theorem proved in~\cite{pure6}.
We need to prove that for all $\sigma''$ with $\sigma''< c-1/\lfloor(\beta-3)/3\rfloor$, there exist $N'', K'', \lambda''$ such that 
$(N'',K'',\sigma'',\lambda'',c)$ forces $H''$.
To prove this, we choose $\sigma'$ with 
$\sigma''< \sigma'<c-1/\lfloor(\beta-3)/3\rfloor$; we use the inductive hypothesis to obtain $K',\lambda'$ such that 
$(N',K',\sigma',\lambda',c)$ forces $H'$; and then we apply the theorem of~\cite{pure6} to deduce what we want. But, crucially, the 
theorem of~\cite{pure6} gives more than this, by exploiting the fact that $u,v$ both have degree one in $H''$. The theorem of~\cite{pure6} gives $N'', K'', \lambda''$ such that 
$(N'',K'',\sigma'',\lambda'',c)$ forces $H''$ in such a way that the vertices $u,v$ of $H''$ appear in the first and last blocks of 
the blockade that contain any vertices of $H''$ (an ``aligned'' copy of $H''$, say).

For the second part, we need to show that for all $\sigma$ with $\sigma< c-1/\lfloor(\beta-3)/3\rfloor$, there exist $N, K, \lambda$ 
such that $(N,K,\sigma,\lambda,c)$ forces $H$.
Choose $\sigma''$ with $\sigma<\sigma''< c-1/\lfloor(\beta-3)/3\rfloor$, and choose $N'', K'', \lambda''$ such that
$(N'',K'',\sigma'',\lambda'',c)$ forces an aligned copy of $H''$. Now we have some blockade 
$\mathcal{A}$ of huge length, and shrinkage at most $\sigma$, and linkage as small as we want.
We know that in every long enough blockade of shrinkage at most $\sigma''$ and linkage at most $\lambda''$
there is an aligned rainbow copy of $H''$.
In particular, there are many aligned $\mathcal{A}$-rainbow copies of $H''$, but we don't know which will be the pairs of blocks containing
the first and last vertices of such a copy.

We need to build a contraption that will allow us to add an $\mathcal{A}$-rainbow handle of the desired length between a pair of blocks of $\mathcal{A}$, when
we don't yet know the right pair of blocks. This we can do, with a device we call a ``bi-grading'', at the cost of increasing the 
shrinkage of the blockade by a small constant (which we can afford, because of the wiggle room between $\sigma, \sigma''$). 
Obtaining this 
device is the main part of the paper, and we omit further details here.

\section{Expansion}

To avoid constantly having to refer
to a blockade $\mathcal{A}$, we define the {\em $\mathcal{A}$-size} of a set $X$ to be $|X|/|A|$, if $X\subseteq A$ for some
block $A$ of $\mathcal{A}$.
For $\gamma,\delta\ge 0$, we say that a blockade $\mathcal{A}=(A_i:i\in I)$ is
{\em $(\gamma,\delta)$-divergent} if there exist distinct $i,j\in I$, and $X\subseteq A_i$ and $Y\subseteq A_j$
such that $|X|\ge \gamma|A_i|$, and $|Y|\ge \delta|A_j|$, and $X$ is anticomplete to $Y$.

Let $\mathcal{A}=(A_i:i\in I)$ be a blockade in $G$. If $J\subseteq I$, we say $(A_i:i\in J)$ is a {\em sub-blockade} of $\mathcal{A}$;
and if $B_i\subseteq A_i$ for each $i\in I$, we say that $(B_i:i\in I)$ is a {\em contraction} of $\mathcal{A}$.
Let $\mathcal{B}=(B_i:i\in J)$ be a contraction of a sub-blockade of $\mathcal{A}$, and let $\kappa$ be 
the minimum of the $\mathcal{A}$-size of $B_i$, over all $i\in J$. We call $\kappa$ the {\em $\mathcal{A}$-size}
of $\mathcal{B}$.

Let us say a blockade $\mathcal{B}=(B_i:i\in I)$ in $G$ is {\em $\tau$-expanding} if 
$$\frac{|N(X)\cap B_j|}{|B_j|}\ge \min\left(\frac{\tau|X|}{|B_i|},\frac14\right)$$
for all distinct
$i,j\in I$, and all $X\subseteq B_i$.

\begin{thm}\label{blockexpand}
Let $G$ be a graph, and let $K\ge 2$ be an integer. Let $\delta>0$ with $\delta K \le 1/4$, and let $\mathcal{A}=(A_i:i\in I)$ be a 
blockade in $G$ of length $K$ that is not $(1/8,\delta)$-divergent.
For each $i\in I$ there exists $B_i\subseteq A_i$ with $|B_i|\ge (1-\delta K)|A_i|$ such that $\mathcal{B}=(B_i:i\in I)$ 
is $(1/(4\delta))$-expanding.
\end{thm}
\Proof
For all $i,j\in I$, let $Z_{i,j}\subseteq A_i$, where $Z_{i,i}=\emptyset$. For all $i\in I$, let $Z_i=\bigcup_{j\in I}Z_{i,j}$, 
and for all distinct $i,j\in I$, let $Y_{i,j}$ denote
the set of vertices in $A_j\setminus Z_j$ that have a neighbour in $Z_{i,j}$.
We say 
such a choice of sets $Z_{i,j}\; (i,j\in I)$ is {\em good} if 
$|Z_{i,j}|/|A_i|< \delta$ and
$|Y_{i,j}|/|A_j|\le |Z_{i,j}|/(3\delta|A_i|)$ for all distinct $i,j\in I$.
\\
\\
(1) {\em If $Z_{i,j}\; (i,j\in I)$ is a good choice of sets, then $|Z_i|\le \delta K|A_i|\le |A_i|/4$ for each $i\in I$, and  $|Y_{i,j}|\le |A_j|/3$
for all distinct $i,j\in I$.}
\\
\\
Since each $|Z_{i,j}|\le \delta|A_i|$,
it follows that $|Z_i|\le \delta K|A_i|\le |A_i|/4$ for each $i\in I$.
Also, for all distinct $i,j\in I$, since $|Y_{i,j}|/|A_j|\le |Z_{i,j}|/(3\delta|A_i|)$ and $|Z_{i,j}|\le \delta|A_i|$, it follows
that $|Y_{i,j}|\le |A_j|/3$. 
This proves (1).

\bigskip

There is a good choice of sets $Z_{i,j}\; (i,j\in I)$, since we may take each $Z_{i,j}=\emptyset$.
Let $Z_{i,j}\; (i,j\in I)$ be a good choice of sets with $\sum_{i,j\in I}|Z_{i,j}|$ maximum (we call this {\em optimality}).
\\
\\
(2) {\em 
Let $i,j\in I$ be distinct, let $X\subseteq A_i\setminus Z_i$, and let $Y$ be the set of vertices in 
$A_j\setminus Z_j$ that have a neighbour in $X$. Then 
$$|Y|/|A_j|\ge \min(|X|/(3\delta|A_i|),1/4)\ge \min(|X|/(4\delta|A_i\setminus Z_i|),1/4).$$
}
To prove the first inequality, we may assume 
that $|Y|/|A_j|< |X|/(3\delta|A_i|)$ (and therefore $X\ne \emptyset$).
Since $|Y_{i,j}|/|A_j|\le |Z_{i,j}|/(3\delta|A_i|)$,
it follows that $|Y\cup Y_{i,j}|/|A_j|\le  |X\cup Z_{i,j}|/(3\delta|A_i|)$.

From optimality,
adding $X$ to $Z_{i,j}$ violates one of the conditions of ``good choice'', and we have just seen that it does not violate 
the second condition. So it violates the first, that is, 
$|Z_{i,j}\cup X|\ge \delta|A_i|$. 
Since 
$\mathcal{A}$ is not $(\delta, 1/8)$-divergent,
fewer than  $|A_j|/8$ vertices in $A_j$ are anticomplete to $Z_{i,j}\cup X$, and so at least $7|A_j|/8$
vertices in $A_j$ have a neighbour in $Z_{i,j}\cup X$. But all such vertices 
belong to $Z_j\cup Y_{i,j}\cup Y$; and since $|Z_j|\le |A_j|/4$ and $|Y_{i,j}|\le |A_j|/3$ by (1),
it follows that $|Y|\ge 7|A_j|/24\ge |A_j|/4$. This proves the first inequality of (2).

Since
$|Z_i|\le |A_i|/4$, and $\delta\le 1/8$ (because $\delta K\le 1/4$ and $K\ge 2$), it follows that 
$$|Y|/|A_j\setminus Z_j|\ge |Y|/|A_j|\ge \min(|X|/(3\delta|A_i|),1/4)\ge \min(|X|/(4\delta|A_i\setminus Z_i|),1/4).$$
This proves (2).

\bigskip

Hence $(A_i\setminus Z_i:i\in I)$ is $(1/(4\delta))$-expanding. 
This proves \ref{blockexpand}.~\bbox

In the next result, we have changed the index set of the blockade from $I$ to $H$, for convenience when we apply it later.
\begin{thm}\label{shortpath}
Let $G$ be a graph, and let $\rho\ge 2$ be an integer. Let $\delta>0$ with $\delta \rho \le 1/4$, and let $\mathcal{A}=(A_i:i\in H)$ be a
blockade in $G$ of length $\rho$, that is not $(\gamma,\delta)$-divergent, where $\gamma\le 1/8$ and $(4\delta)^{\rho}|G|\le 3$.
For each $i\in H$ let $B_i\subseteq A_i$ be as in \ref{blockexpand}, and
$\mathcal{B}=(B_i:i\in H)$.
Let $i_1,i_2\in H$ be distinct, let $v\in B_{i_1}$, and let
$Y\subseteq B_{i_2}$ with $|Y|\ge \gamma|A_{i_2}|$.
Then there is a $\mathcal{B}$-rainbow path with one end $v$ and the other end in $Y$.
\end{thm}
\Proof By hypothesis, $|A_i\setminus B_i|\le \delta \rho|A_i|\le |A_i|/4$ for each $i\in H$, and $\mathcal{B}$ is
$\tau$-expanding, where $\tau=1/(4\delta)$.
Without loss of generality we may assume that $H=\{1\LL \rho\}$, where $i_1=1$ and $i_2=\rho$.
For $1\le i\le \rho-1$, let $X_i$ be the set of vertices in $B_{i}$ that can be joined to $v$
by a $(B_1\LL B_i)$-rainbow path. Since for $i\ge 2$,
$X_i$ contains all vertices
in $B_{i}$ that have a neighbour in $X_{i-1}$, it follows that
$$\frac{|X_i|}{|B_i|}\ge \min\left ( \frac{\tau|X_{i-1}|}{|B_{i-1}|}, \frac14\right),$$
and since $\tau\ge 1$, it follows that $|X_i|/|B_i|\ge \min\left (\tau^{i-1}/|B_{1}|, 1/4\right).$
In particular, since
$$\frac{\tau^{\rho-1}}{|B_1|}\ge \frac{\tau^{\rho-1}}{|G|}\ge \frac{4\delta}{3},$$
(because $(4\delta)^{\rho}|G|\le 3$) it follows that
$$\frac{|X_{\rho-1}|}{|B_{\rho-1}|}\ge \min \left(\frac{4\delta}{3}, \frac{1}{4}\right)=\frac{4\delta}{3};$$
and so $|X_{\rho-1}|/|A_{\rho-1}|\ge \delta$, since
$|B_{\rho-1}|\ge 3|A_{\rho-1}|/4$. Since $\mathcal{A}$ is not $(\gamma,\delta)$-divergent, there are fewer than $\gamma|A_{\rho}|$
vertices in $A_{\rho}$ that have no neighbour in $X_{\rho-1}$; and so $|B_{\rho}\setminus X_{\rho}|< \gamma|A_{\rho}|$.
Since $Y\subseteq B_{\rho}$ and $|Y|\ge \gamma|A_{\rho}|$, it follows that $X_{\rho}\cap Y\ne \emptyset$.
This proves \ref{shortpath}.~\bbox

A {\em levelling} in a graph $G$ is a sequence $(L_0\LL L_k)$ of pairwise disjoint subsets of $V(G)$, where $k\ge 1$, 
with the following properties:
\begin{itemize}
\item $|L_0|=1$;
\item for $1\le i\le k$, $L_{i-1}$ covers $L_i$; and 
\item for $2\le i\le k$, $L_0\cupcup L_{i-2}$ is anticomplete to $L_i$.
\end{itemize}
We call the unique member of $L_0$ the {\em apex} of the levelling, and $L_k$ is its {\em base}, and $k$ is its {\em height}.

Let $\mathcal{L}=(L_0\LL L_k)$ be a levelling in $G$, and let $C\subseteq V(G)$. We say that $\mathcal{L}$ {\em reaches} $C$
if $(L_0,L_1\LL L_k, C)$ is a levelling.

Now let $\mathcal{B}=(B_i:i\in I)$ be a blockade in $G$, and let $\mathcal{L}=(L_0\LL L_k)$ be a levelling in $G$. We say that 
$\mathcal{L}$ is {\em $\mathcal{B}$-rainbow} if for $0\le i\le k$, there exists $h_i\in I$ such that
$h_0\LL h_{k}$ are all distinct, and
$L_i\subseteq B_{h_i}$ for $0\le i\le k$.

\begin{thm}\label{getlevelling}
Let $\tau\ge 6$, and let $\mathcal{B}=(B_i:i\in I)$ be a blockade in a graph $G$. Let $h_0\in I$, and let
$(B_i:i\in I\setminus \{h_0\})$ be $\tau$-expanding.
Let $\rho$ be an integer such that $(\tau/2)^{\rho-1}\ge |G|$, let $H\subseteq I$ with $h_0\in H$ and $|H|=\rho$, and
let $v\in B_{h_0}$, with a neighbour in $\bigcup_{h\in H\setminus \{h_0\}}B_h$. Then there exist $J\subseteq I\setminus H$
with $|J|\ge |I|/\rho-1$, and a
levelling $\mathcal{L}=(L_0\LL L_{k-1}, L_k)$ for some $k\le \rho$, with apex $v$, such that $(L_0\LL L_{k-1})$ is
$(B_i:i\in H)$-rainbow, and
$|L_k\cap B_j|\ge |B_j|/(4\rho)$ for all $j\in J$.
\end{thm}
\Proof
Let $L_0:=\{v\}$. Since $v$ has  a neighbour in $\bigcup_{h\in H\setminus \{h_0\}}B_h$, we may choose $h\in H\setminus \{h_0\}$ such that $|N(v)\cap B_{h}|/|B_h|$
is maximum, and we set $h_1=h$ and $m_1=|N(v)\cap B_{h}|/|B_h|$, and $L_1=N(v)\cap B_{h}$. Thus $m_1\ge 1/|G|$.
We define $t$ and $h_2\LL h_{t}$ inductively as follows.
Assume inductively that for some $i\le \rho-1$, we have already defined
$h_0\LL h_i$, and $L_{0}\LL L_{i}$, and $m_1\LL m_i$,
with the properties that:
\begin{itemize}
\item $h_0,h_1\LL h_i\in H$ are all distinct;
\item $L_g\subseteq B_{h_g}$ for $0\le g\le i$;
\item $(L_0\LL L_i)$ is a levelling;
\item $m_{g}=|N(L_{g-1})\cap B_{h_{g}}|/|B_{h_{g}}|$ for $1\le g\le i$;
\item $|N(L_{g-1})\cap B_h|/|B_h|\le m_{g}$ for all $g$ with $1\le g\le i$ and all $h\in H\setminus \{h_0,h_1\LL h_i\}$;
\item $m_g\ge (\tau-3)(m_1+\cdots+ m_{g-1})$ and $|L_g|\ge \left(1-\frac{2}{\tau}\right)m_g|B_{h_g}|$ for $1\le g\le i$.
\end{itemize}

If $|L_i|\ge |B_{h_i}|/(4\tau)$, let $t:=i$ and the inductive definition is complete. Otherwise we proceed as follows.
Since $|L_i|< |B_{h_i}|/(4\tau)$, and $|L_i|\ge (1-2/\tau)m_i|B_{h_i}|$, it follows that
$(1-2/\tau)m_i<1/(4\tau)$, and so $m_i<1/(4(\tau-2))$.
But
$$m_g\ge (\tau-3)(m_1+\cdots+ m_{g-1})\ge (\tau-3)m_{g-1}$$
for $2\le g\le i$, and consequently $m_i\ge (\tau-3)^{i-1}m_1\ge (\tau-3)^{i-1}/|G|$; and therefore
$$|G|>4(\tau-2)(\tau-3)^{i-1}\ge (\tau/2)^{i}.$$ 
Since $|G|\le (\tau/2)^{\rho-1}$, it follows that
$i<\rho-1$, and so
$$|\{h_0,h_1\LL h_i\}|<\rho=|H|.$$

Choose $h\in H\setminus \{h_0,h_1\LL h_i\}$ with $|N(L_i)\cap B_{h}|/|B_h|$ maximum, and define $h_{i+1}=h$
and $m_{i+1}=|N(L_i)\cap B_{h}|/|B_h|$. Let $L_{i+1}$ be the set of vertices in $N(L_i)\cap B_{h_{i+1}}$ that have no
neighbour in $L_0\cupcup L_{i-1}$. Thus the first five conditions above are satisfied. It remains only to check the final condition,
that is, $m_{i+1}\ge (\tau-3)(m_1+\cdots+ m_{i})$ and
$|L_{i+1}|\ge (1-2/\tau)m_{i+1}|B_{h_{i+1}}|$.

Since $(B_i:i\in I\setminus \{h_0\})$ is $\tau$-expanding and $\frac{|L_i|}{|B_{h_i}|}< \frac{1}{4\tau}$, it follows that
$m_{i+1}\ge \tau\frac{|L_i|}{|B_{h_i}|}$.
But $|L_i|\ge \left(1-\frac{2}{\tau}\right)m_i|B_{h_i}|$, and so 
$$m_{i+1}\ge (\tau-2)m_i= (\tau-3)m_i +m_i\ge (\tau-3)m_i+ (\tau-3)(m_1+\cdots+ m_{i-1}).$$
This proves the first part of the final condition. For the second part, we observe
that for $0\le g\le i-1$, the number of vertices in $N(L_i)\cap B_{h_{i+1}}$ that have a neighbour in $L_g$ is at most
$m_g|B_{h_{i+1}}|$, because of the choice of $h_{g+1}$. Thus
$$\frac{|L_{i+1}|}{|B_{h_{i+1}}|}\ge m_{i+1}-(m_1+\cdots+m_i)\ge \left(1-\frac{1}{\tau-3}\right)m_{i+1}\ge \left(1-\frac{2}{\tau}\right)m_i.$$
Thus the final condition holds. This completes the inductive definition. We see that $t\le \rho-1$, and 
$|L_t|\ge |B_{h_t}|/(4\tau)$.

For $0\le i\le t$, let $n_i$ be the number of $j\in I\setminus H$ such that at least $(i+1)|B_j|/(4(t+1))$ vertices in $B_j$ have a neighbour
in $L_0\cupcup L_i$.
We see that $|L_{t}|\ge |B_{h_{t}}|/(4\tau)$; and so, for each $j\in I\setminus H$, since $(B_i:i\in I\setminus \{h_0\})$ is $\tau$-expanding,
it follows that at least $|B_j|/4$ vertices in $B_j$ have a neighbour in $L_{t}$; and so $n_{t}=|I|-|H|$. Choose $k\in \{1\LL t+1\}$
minimum such that $n_{k-1}\ge (k/(t+1))(|I|-|H|)$. It follows from the minimality of $k$ that there are
at least $(|I|-|H|)/(t+1)$ values of $j\in I\setminus H$ such that at least $k|B_j|/(4(t+1))$ vertices in $B_j$ have a neighbour
in $L_0\cupcup L_{k-1}$, and at most $(k-1)|B_j|/(4(t+1))$ vertices in $B_j$ have a neighbour
in $L_0\cupcup L_{k-2}$ (this last statement is vacuously true if $k=1$). Let $J$ be the set of all such values of $j$;
thus $|J|\ge (|I|-|H|)/(t+1)\ge (|I|-|H|)/\rho$. For each $j\in J$,
let $C_j\subseteq B_j$ be the set of all vertices
in $B_j$ that have a neighbour in $L_0\cupcup L_{k-1}$ and have no neighbour in $L_0\cupcup L_{k-2}$ (and therefore
have a neighbour in  $L_{k-1}$). It follows that $|C_j|\ge |B_j|/(4(t+1))\ge |B_j|/(4\rho)$ for each $j\in J$.
Let $L_{k}:=\bigcup_{j\in J}C_j$; then $(L_0\LL L_{k})$ is a levelling with the properties required.
This proves \ref{getlevelling}.~\bbox

\section{Gradings}

Let $\mathcal{B}=(B_i:i\in I)$ be a blockade in a graph $G$, and let $\mathcal{L}=(L_0\LL L_k)$ be a 
levelling. We say that $\mathcal{L}$ {\em grades} $\mathcal{B}$ if $\mathcal{L}$ reaches $\bigcup_{i\in I}B_i$, and 
$I$ can be written as $\{i_1\LL i_n\}$ (not necessarily listed in increasing order), such that 
for $1\le g\le n$ there exists 
$Y\subseteq L_{k}$ that
covers $\bigcup_{g\le h\le n}B_{i_h}$ and
is anticomplete to $\bigcup_{1\le h<g}B_{i_h}$.

\begin{thm}\label{getgrading}
Let $\tau\ge 6$, and let $\mathcal{B}=(B_i:i\in I)$ be a $\tau$-expanding blockade with linkage at most $1/(8|I|)$ 
in a graph $G$.
Let $\rho$ be an integer such that $(\tau/2)^{\rho-1}\ge |G|$, let $H\subseteq I$ with $|H|=\rho$, let $h_0\in H$, and 
let $v\in B_{h_0}$. Then there exist $J\subseteq I\setminus H$ with $|J|= \lceil |I|/\rho\rceil -1$, and 
$C_j\subseteq B_j$ with $|C_j|\ge |B_j|/(8|I|)$ for all $j\in J$,
and a $(B_i:i\in H)$-rainbow
levelling $\mathcal{L}$ with apex $v$ that
grades $(C_j:j\in J)$.
\end{thm}
\Proof
By \ref{getlevelling}, there exists $J\subseteq I\setminus H$ with $|J|\ge |I|/\rho-1$, and a 
levelling $\mathcal{L}=(L_0\LL L_k)$ for some $k\le \rho$, with apex $v$, such that $(L_0\LL L_{k-1})$ is $(B_i:i\in H)$-rainbow, and
$|L_k\cap B_j|\ge |B_j|/(4\rho)$ for all $j\in J$.  We may choose $J$ with $|J|= \lceil |I|/\rho\rceil -1$.

Let $|J|=n$, and $Y_0=\emptyset$. We define $Y_1\LL Y_n\subseteq L_{k-1}$, and distinct $j_1\LL j_n\in J$ inductively as follows.
Let $1\le i\le n$, and suppose that $Y_1\LL Y_{i-1}$ and $j_1\LL j_{i-1}$ have been defined, with $Y_0\subseteq Y_1\subseteq \cdots\subseteq Y_{i-1}$.
Choose $Y_{i}\subseteq L_{k-1}$ including $Y_{i-1}$, minimal such that $|N(Y_{i})\cap L_k\cap B_j|\ge i|B_j|/(4\rho n)$
for some $j\in J\setminus \{j_1\LL j_{i-1}\}$; and let $j_{i}:=j$ for some such choice of $j$. (This is possible since $i\le n$,
and $|L_k\cap B_j|\ge |B_j|/(4\rho)$, and $L_{k-1}$ covers $L_k\cap B_j$.) This completes the inductive definition.
\\
\\
(1) {\em For $1\le i\le n$, $Y_i\ne Y_{i-1}$, and 
$$|N(Y_{i})\cap L_k\cap B_{j}|\le \left(\frac{1}{8|I|}+\frac{i}{4\rho n}\right) |B_{j}|$$ 
for all $j\in J\setminus \{j_1\LL j_{i-1}\}$.}
\\
\\
We prove both statements simultaneously by induction on $i$. Thus, we assume both statements hold for $i-1$ (if $i>1$).
To prove the first statement holds for $i$, we may assume that $i\ge 2$, because clearly $Y_1\ne Y_0$.
Since 
$$|N(Y_{i-1})\cap L_k\cap B_{j}|\le \left(\frac{1}{8|I|}+\frac{i-1}{4\rho n}\right) |B_j|$$
for all $j\in J\setminus \{j_1\LL j_{i-2}\}$, and in particular for $j=j_i$, and since
$$|N(Y_{i})\cap L_k\cap B_{j_i}|\ge \frac{i}{4\rho n}|B_{j_i}|> \left(\frac{1}{8|I|}+\frac{i-1}{4\rho n}\right) |B_{j_i}|$$
(because $n\le |I|/\rho$)
it follows that $Y_i\ne Y_{i-1}$. This proves the first statement of (1). To prove the second statement, since $Y_i\ne Y_{i-1}$, we may choose $u\in Y_i\setminus Y_{i-1}$;
and by the minimality of $Y_i$, for each $j\in J\setminus \{j_1\LL j_{i-1}\}$ it follows that 
$$|N(Y_{i}\setminus \{u\})\cap L_k\cap B_{j}|< \frac{i}{4\rho n}|B_{j}|.$$
But $\mathcal{B}$ has linkage at most $1/(8|I|)$, and so at most $|B_{j}|/(8|I|)$ vertices in $B_{j}$ are adjacent to $u$; and 
therefore
$$|N(Y_{i})\cap L_k\cap B_{j}|\le \frac{i}{4\rho n}|B_{j}|+\frac{1}{8|I|}|B_{j}|.$$
This proves the second statement and so proves (1).

\bigskip

For $1\le i\le n$, let $C_{j_i}$ be the set of vertices in $B_{j_i}\cap L_k$ that have a neighbour in $Y_i$ and have no neighbour in $Y_{i-1}$.
Thus by (1), 
$$\frac{|C_i|}{|B_{j_i}|}\ge \frac{i}{4\rho n}- \left(\frac{1}{8|I|}+\frac{i-1}{4\rho n}\right) \ge \frac{1}{4\rho n}- \frac{1}{8|I|}\ge \frac{1}{8\rho n}\ge \frac{1}{8|I|}.$$
For $1\le i\le n$,
$Y_i$ covers $C_{j_1}\LL C_{j_i}$, and is anticomplete to $C_{j_{i+1}}\LL C_{j_n}$.
Hence  $(L_0\LL L_{k-1})$ grades $(C_j:j\in J)$.
This proves \ref{getgrading}.~\bbox

Now we come to the main result of this section, an extension of \ref{getgrading}. We say two levellings $(L_0\LL L_k)$
and $(M_0\LL M_m)$ in a graph $G$ are {\em parallel} if 
$L_0=M_0$, and
$M_1\cupcup M_{m}$ is disjoint from and anticomplete to $L_1\cupcup L_{k}$.
Let $\mathcal{B}=(B_i:i\in I)$ be a blockade in a graph $G$, and let $\mathcal{L}=(L_0\LL L_k)$ be a
levelling. We say that $\mathcal{L}$ grades $\mathcal{B}$ {\em forwards} if $\mathcal{L}$ reaches $\bigcup_{i\in I}B_i$, and
for each $j\in I$ there exists
$Y\subseteq L_{k}$ that
covers $\bigcup_{i\in I, i\ge j}B_i$ and
is anticomplete to $\bigcup_{i\in I, i<j}B_i$.

Let $\mathcal{C}=(C_j:j\in J)$ be a blockade in a graph $G$. Let $\mathcal{L}=(L_0\LL L_k)$ and $\mathcal{M}=(M_0\LL M_m)$ be 
parallel levellings both reaching 
$\bigcup_{j\in J}C_j$. Let $\mathcal{L}$ grade $\mathcal{C}$ forwards. In these circumstances we say that
$(\mathcal{L}, \mathcal{M},\mathcal{C})$ is a {\em bi-levelling}, and the blocks of $\mathcal{C}$ are called the {\em blocks} 
of the bi-levelling. We call the sum of the heights of $\mathcal{L}$ and $\mathcal{M}$
the {\em height} of the bi-levelling, and its {\em length} is the
length of $\mathcal{C}$.

Let $\mathcal{A}=(A_i:i\in I)$ be another
blockade in $G$. 
We say that the bi-levelling $(\mathcal{L}, \mathcal{M},\mathcal{C})$ is {\em $\mathcal{A}$-rainbow} if
\begin{itemize}
\item each of the sets $C_j\;(j\in J)$, $L_j\;(0\le j\le k)$, $M_j\;(1\le j\le m)$ is a subset of one of the sets $A_i\;(i\in I)$; and
\item for each $i\in I$, $A_i$ includes at most one of the sets $C_j\;(j\in J)$, $L_j\;(0\le j\le k)$, $M_j\;(1\le j\le m)$.
\end{itemize}
We stress that here the order of the blocks of $\mathcal{A}$ is immaterial. In particular, the orders of the blocks of 
$\mathcal{C}$ may be different from the order of the corresponding blocks of $\mathcal{A}$.

The following is immediate but crucial, so we state it explicitly.
\begin{thm}\label{bilevelpath}
Let $(\mathcal{L}, \mathcal{M},\mathcal{C})$ be an $\mathcal{A}$-rainbow bi-levelling in $G$ of height $\ell$, and let $x,y$ belong 
to the bases of 
$\mathcal{M},\mathcal{L}$ respectively. There is an induced path $P$ between $x,y$, of length $\ell$, such that each of its vertices 
belongs to a different block of $\mathcal{A}$, and each of these blocks of $\mathcal{A}$ includes no block of $\mathcal{C}$; 
and no internal vertex of $P$ has a neighbour in any block of $\mathcal{C}$.
\end{thm}

If $\mathcal{A}=(A_i:i\in I)$ is a blockade, and $\phi:J\rightarrow I$ is a bijection, where $J$ is a set of integers, let 
$B_j:=A_{\phi(j)}$ for each $j\in J$; then $\mathcal{B}=(B_j:j\in J)$ is also a blockade, with the same blocks, but in a different order.
We say that $\mathcal{B}$ is obtained from $\mathcal{A}$ by {\em re-indexing}. Some of our results are invariant under 
re-indexing the blockade in question, and it can simplify notation to take advantage of this.

\begin{thm}\label{getbilevel}
Let $k\ge 1$ and $\rho\ge 2$ be integers. Let $K:=k\rho^4$, let $\gamma,\lambda>0$ with $\lambda\le 1/(512\rho^2K)$ and $\gamma\le 3/(256K)$,
let $G$ be a graph, let $\delta>0$ with 
$\delta\le 3\rho/(128 K^2)$ and
$(256K\delta/3)^{\rho-1}|G| \le 1$, and let $\mathcal{A}=(A_i:i\in I)$ be a
blockade in $G$ of length $K$, with linkage at most $\lambda$, that is not $(\gamma,\delta)$-divergent.
Then there is an $\mathcal{A}$-rainbow bi-levelling $(\mathcal{L}, \mathcal{M},\mathcal{C})$ with length $k$ and height
at most $3\rho-3$, such that the $\mathcal{A}$-size of $\mathcal{C}$ is at least $1/(64\rho^3K)$.
\end{thm}

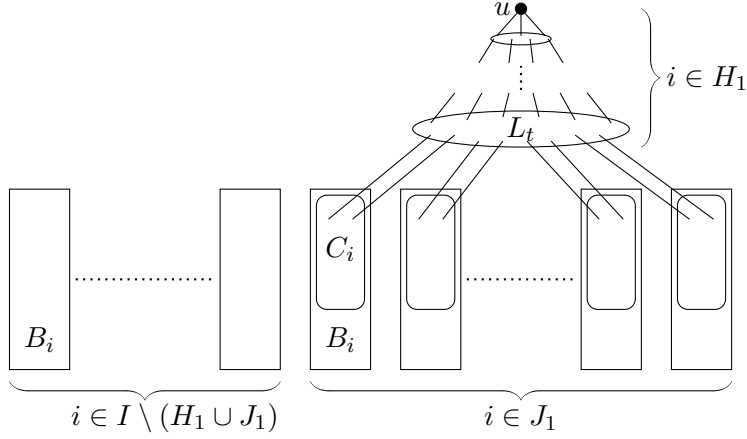
\begin{figure}[h!]
\centering

\begin{tikzpicture}[scale=0.8,auto=left]
\draw[] (-10,0) rectangle (-9,3);
\draw[] (-8.5,0) rectangle (-7.5,3);
\draw[] (-5.5,0) rectangle (-4.5,3);
\draw[] (-4,0) rectangle (-3,3);

\draw[] (-15,0) rectangle (-14,3);
\draw[] (-11.5,0) rectangle (-10.5,3);

\def\r{.1}
\draw[rounded corners] (-10+\r,1) rectangle (-9-\r,3-\r);
\draw[rounded corners] (-8.5+\r,1) rectangle (-7.5-\r,3-\r);
\draw[rounded corners] (-5.5+\r,1) rectangle (-4.5-\r,3-\r);
\draw[rounded corners] (-4+\r,1) rectangle (-3-\r,3-\r);

\draw (-6.5,4) ellipse (1.8 and .3);
\draw (-6.5,5.5) ellipse (.5 and .1 );

\tikzstyle{every node}=[inner sep=1.5pt, fill=black,circle,draw]
\node (u) at (-6.5,6) {};

\draw [decorate,
        decoration = {calligraphic brace, amplitude=7pt}] (-3,-.2) --  (-10,-.2);
\draw [decorate,
        decoration = {calligraphic brace, amplitude=7pt}] (-10.5,-.2) --  (-15,-.2);
\draw [decorate,
        decoration = {calligraphic brace, amplitude=7pt}] (-4.5,6) --  (-4.5, 3.7);
\tikzstyle{every node}=[]
\draw (-6.5,-.8) node []           {$i\in J_1$};
\draw (-12.25,-.8) node []           {$i\in I\setminus (H_1\cup J_1)$};
\draw (-3.4,4.85) node []           {$i\in H_1$};
\draw (-9.5,.5) node []           {$B_i$};
\draw (-9.5,2) node []           {$C_i$};
\draw (-14.5,.5) node []           {$B_i$};
\draw (u) node [left]           {$u$};
\draw (-6.5,4) node []       {$L_t$};

\draw (u) to (-6.5,5.55);
\draw (u) to (-6.9,5.55);
\draw (u) to (-6.1,5.55);
\draw[]  (-6.95,5.5) -- (-7.2,5.1);
\draw[]  (-6.65,5.5) -- (-6.7,5.1);
\draw[]  (-6.35,5.5) -- (-6.3,5.1);
\draw[]  (-6.05,5.5) -- (-5.8,5.1);

\draw[]  (-7.6,4.6) -- (-8,4.1);
\draw[]  (-7.15,4.6) -- (-7.4,4.15);
\draw[]  (-6.7,4.6) -- (-6.8,4.2);
\draw[]  (-6.3,4.6) -- (-6.2,4.2);
\draw[]  (-5.85,4.6) -- (-5.6,4.15);
\draw[]  (-5.4,4.6) -- (-5,4.1);

\draw[dotted, thick] (-6.5,4.65) to (-6.5,5.05);
\draw[dotted, thick] (-7.4,1.5) to (-5.6, 1.5);
\draw[dotted, thick] (-13.9,1.5) to (-11.6,1.5);
\draw (-8, 3.9) to (-9.7,2.5);
\draw (-7.6, 3.9) to (-9.3,2.5);
\draw (-7.2, 3.8) to (-8.2,2.5);
\draw (-6.8, 3.8) to (-7.8,2.5);
\draw (-6.4, 3.8) to (-5.2,2.5);
\draw (-6, 3.8) to (-4.8,2.5);
\draw (-5.6, 3.9) to (-3.7,2.5);
\draw (-5.2, 3.9) to (-3.3,2.5);

\end{tikzpicture}
\caption{$\mathcal{B}$, $\mathcal{C}$, and $\mathcal{L}$. Everything is $\mathcal{A}$-rainbow, and $\mathcal{L}$ grades $\mathcal{C}$ forwards (not shown). }\label{fig:getbilevel1}
\end{figure}

\Proof
By \ref{blockexpand}, for each $i\in I$ there exists $B_i\subseteq A_i$ with $|A_i\setminus B_i|\le \delta K|A_i|$, such that
$\mathcal{B}=(B_i:i\in I)$ is $(1/(4\delta))$-expanding. 
It follows that $\mathcal{B}$ has linkage at most $4\lambda/3$, since
$|B_i|\ge (1-\delta K)|A_i|\ge 3|A_i|/4$ for each $i\in I$.
Choose $H_1\subseteq I$ with $|H_1|=\rho$. Since $8(1-\delta K)^{-1}\gamma |I|\le 1$,
by \ref{getgrading} applied to $\mathcal{B}$, taking $\tau=1/(4\delta)$, there exist $J_1\subseteq I\setminus H_1$ 
with $|J_1|= \lceil |I|/\rho\rceil -1$,
and $C_j\subseteq B_j$ with $|C_j|\ge |B_j|/(8|I|)$ for all $j\in J_1$,
and a $(B_i:i\in H_1)$-rainbow
levelling $\mathcal{L}=(L_0\LL L_t)$ that
grades $\mathcal{C}=(C_j:j\in J_1)$. (See figure \ref{fig:getbilevel1}.) Thus $\mathcal{L}$ has height at most $\rho-1$.
The statement of the theorem is invariant under re-indexing the blocks of $\mathcal{A}$; and so we may assume without loss of
generality that $\mathcal{L}$ grades $\mathcal{C}$ forwards (by re-indexing appropriately $\mathcal{C}$, and
correspondingly $\mathcal{B},\mathcal{A}$).

Since $\mathcal{A}$ is not $(\gamma,\delta)$-divergent, and 
$$|C_j|\ge \frac{1}{8|I|}|B_j|\ge \frac{3}{32|I|}|A_j|$$ 
for each $j\in J_1$, it follows that
$\mathcal{C}$ is not $(32|I|\gamma/3,32|I|\delta/3)$-divergent. Since $32|I|\gamma/3\le 1/8$, and $32|I|\delta |J_1|/3\le 1/4$, 
it follows from \ref{blockexpand}  applied to $(C_j:j\in J_1)$ (with $\delta$ replaced by $32|I|\delta/3$) that for each $j\in J_1$ 
there exists $D_j\subseteq C_j$ with $|D_j|\ge 3|C_j|/4$, such that $\mathcal{D}=(D_j:j\in J_1)$ is $(3/(128|I|\delta))$-expanding.
Thus 
$$|D_j|\ge \frac{3}{4}|C_j|\ge \frac{3}{32|I|} |B_j|\ge \frac{9}{128|I|}|A_j|.$$

Let the apex of $\mathcal{L}$, $u$ say, belong to $B_{h_1}$ where $h_1\in H_1$.
Choose $H_2\subseteq I\setminus (H_1\cup J_1)$ with cardinality 
$\rho-2$.
(Since $|H_1|\le \rho$, and $|J_1|\le |I|/\rho$, it follows that 
$$|I\setminus (H_1\cup J_1)|\ge |I|(1-1/\rho)-\rho\ge \rho-2,$$
so this is possible.) Let $H_3\subseteq J_1$ be the set of the $\rho-1$ smallest members of $J_1$, and let $j_1\in H_3$.
Since
$$|D_{j_1}|\ge \frac{9}{128|I|}|A_{j_1}|\ge \gamma |A_{j_1}|,$$
it follows from \ref{shortpath} applied to $(A_j:j\in H_2\cup \{h_1,j_1\})$ that there is a $(B_i:i\in H_2\cup \{h_1,j_1\})$-rainbow path
with one end $u$ and the other in $D_{j_1}$.
Its length is at most $\rho-1$. Consequently there is a $(B_i:i\in H_2\cup \{h_1\})$-rainbow
path $P$ of minimum length such that 
one end is $u$, and the other end, $v$ say, has a neighbour in $\bigcup_{j\in H_3}D_{j}$.
It follows that $P$ has length at most $\rho-2$, and 
$V(P)\setminus \{v\}$ is anticomplete to 
$\bigcup_{j\in H_3}D_{j}$.

\begin{figure}[h!]
\centering

\begin{tikzpicture}[scale=0.8,auto=left]
\draw[] (-18,0) rectangle (-17,3);
\draw[] (-15.5,0) rectangle (-14.5,3);
\draw[] (-14,0) rectangle (-13,3);
\draw[] (-11.5,0) rectangle (-10.5,3);

\draw[] (-10,0) rectangle (-9,3);
\draw[] (-7.5,0) rectangle (-6.5,3);
\draw[] (-6,0) rectangle (-5,3);
\draw[] (-3.5,0) rectangle (-2.5,3);

\draw[] (-2,0) rectangle (-1,3);
\draw[] (.5,0) rectangle (1.5,3);

\def\s{.25}

\def\r{.1}
\draw[rounded corners] (-10+\r,.5) rectangle (-9-\r,3-\r);
\draw[rounded corners] (-7.5+\r,.5) rectangle (-6.5-\r,3-\r);
\draw[rounded corners] (-6+\r,.5) rectangle (-5-\r,3-\r);
\draw[rounded corners] (-3.5+\r,.5) rectangle (-2.5-\r,3-\r);

\def\r{.2}
\draw[rounded corners] (-10+\r,1) rectangle (-9-\r,3-\r);
\draw[rounded corners] (-7.5+\r,1) rectangle (-6.5-\r,3-\r);
\draw[rounded corners] (-6+\r,1) rectangle (-5-\r,3-\r);
\draw[rounded corners] (-3.5+\r,1) rectangle (-2.5-\r,3-\r);

\def\r{.3}
\draw[rounded corners] (-6+\r,1.5) rectangle (-5-\r,3-\r);
\draw[rounded corners] (-3.5+\r,1.5) rectangle (-2.5-\r,3-\r);

\draw (-6.5+\s,4) ellipse (1.8 and .3);
\draw (-6.5+\s,5) ellipse (.5 and .1 );
\draw (-9.5,2) ellipse (.1 and .4 );
\draw (-7,2) ellipse (.3 and .6 );

\tikzstyle{every node}=[inner sep=1.5pt, fill=black,circle,draw]
\node (u) at (-6.5+\s,5.5) {};
\node (w) at (-13.5,2) {};
\node (x) at (-12.8,2) {};
\node (y) at (-11.7,2) {};
\node (v) at (-11,2) {};

\draw[] (u) to [bend right = 20] (w);

\draw[] (w)--(x);
\draw[dashed, thick] (x)--(y);
\draw[] (y)--(v);

\draw [decorate,
        decoration = {calligraphic brace, amplitude=7pt}] (-6.5,-.2) --  (-10,-.2);
\draw [decorate,
        decoration = {calligraphic brace, amplitude=7pt}] (-2.5,-.2) --  (-6,-.2);
\draw [decorate,
        decoration = {calligraphic brace, amplitude=7pt}] (1.5,-.2) --  (-2,-.2);

\draw [decorate,
        decoration = {calligraphic brace, amplitude=7pt}] (-4.5+\s,5.5) --  (-4.5+\s, 3.7);
\draw [decorate,
        decoration = {calligraphic brace, amplitude=7pt}] (-10.5,-.2) --  (-14,-.2);
\draw [decorate,
        decoration = {calligraphic brace, amplitude=7pt}] (-14.5,-.2) --  (-18,-.2);
\tikzstyle{every node}=[]
\draw (-4.25,-.8) node []           {$i\in J_2$};
\draw (-8.25,-.8) node []           {$i\in H_3$};
\draw (-12.25,-.8) node []           {$i\in H_2$};
\draw (-.25,-.8) node []           {$i\in J_1\setminus(H_3\cup J_2)$};
\draw (-3.4+\s,4.6) node []           {$i\in H_1$};
\draw (-16.25,-.8) node []           {$i\in I\setminus (H_1\cup H_2\cup J_1)$};
\draw (-9.5,.25) node []           {$B_i$};
\draw (-9.5,.75) node []           {$C_i$};
\draw (-9.5,1.25) node []           {$D_i$};
\draw (-13.5,.5) node []           {$B_i$};
\draw (u) node [right]           {$u$};
\draw (v) node [above]           {$v$};
\draw (-6.5+\s,4) node []       {$L_t$};
\draw (-7,2) node []       {\scriptsize $M_m$};
\draw (-5.5,2.25) node []       {\scriptsize $E_i$};
\draw (-12.25,2.25) node []           {$P$};
\draw (u) to (-6.5+\s,5);
\draw (u) to (-6.9+\s,5);
\draw (u) to (-6.1+\s,5);
\draw[dotted, thick] (-6.5+\s,4.3) to (-6.5+\s,4.8);
\draw[dotted, thick] (-8.7,2) to (-8, 2);
\draw[dotted, thick] (-4.9,1.5) to (-3.6, 1.5);
\draw[dotted, thick] (-.9,1.5) to (.4, 1.5);
\draw[dotted, thick] (-16.9,1.5) to (-15.6,1.5);
\draw[dotted, thick] (-12.9,1.5) to (-11.6,1.5);

\draw (-8+\s, 3.9) to (-9.7,2.5);
\draw (-7.6+\s, 3.9) to (-9.3,2.5);
\draw (-7.2+\s, 3.8) to (-7.2,2.5);
\draw (-6.8+\s, 3.8) to (-6.8,2.5);
\draw (-6.4+\s, 3.8) to (-5.6,2.5);
\draw (-6+\s, 3.8) to (-5.4,2.5);
\draw (-5.6+\s, 3.9) to (-3.1,2.5);
\draw (-5.2+\s, 3.9) to (-2.9,2.5);

\draw (v) to (-9.55, 2.3);
\draw (v) to (-9.55, 2);
\draw (v) to (-9.55, 1.7);
\draw (-9.45, 2.3) to (-8.8,2.4);
\draw (-9.45, 2.1) to (-8.8,2.15);
\draw (-9.45, 1.9) to (-8.8,1.85);
\draw (-9.45, 1.7) to (-8.8,1.6);
\draw (-7.85, 2.4) to (-7,2.5);
\draw (-7.85, 2.2) to (-7.1,2.25);
\draw (-7.85, 2) to (-7.2,2);
\draw (-7.85, 1.8) to (-7.1,1.75);
\draw (-7.85, 1.6) to (-7,1.5);

\end{tikzpicture}
\caption{$P$, $\mathcal{D}$, $\mathcal{E}$ and $\mathcal{M}$. $\mathcal{M}$ reaches $\mathcal{E}$ (not shown). }\label{fig:getbilevel2}
\end{figure}
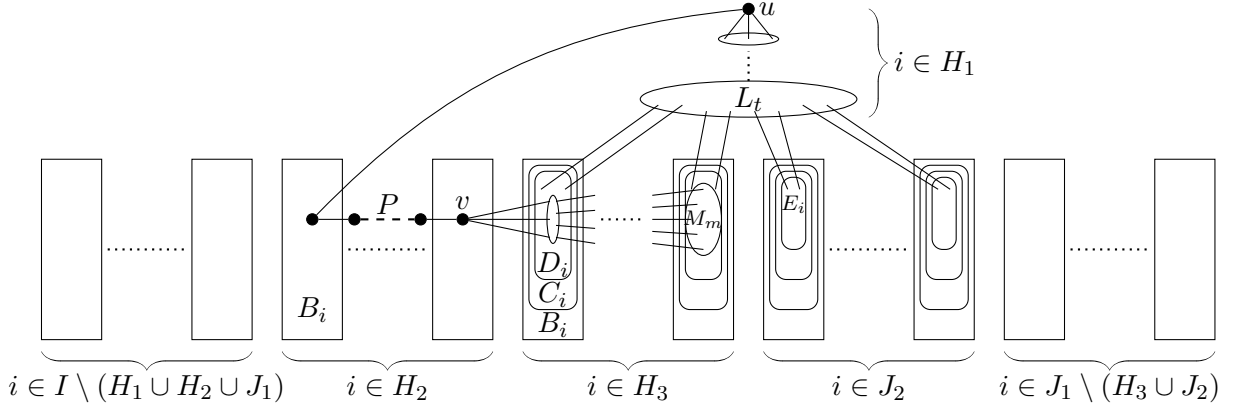

Since $\mathcal{L}$ grades $\mathcal{D}$ and hence 
$\mathcal{L}$ reaches $\bigcup_{i\in J_1}D_i$, it follows that $u$ has no neighbour in  
$\bigcup_{j\in H_3}D_{j}$, and so $v\ne u$. Let $v\in B_{h_2}$ say, where $h_2\in H_2$.
Now $\mathcal{D}$ is $(3/(128|I|\delta))$-expanding, and $3/(128|I|\delta)\ge 6$, and $(3/(256|I|\delta))^{\rho-1}\ge |G|$.
Define $D_{h_2}=B_{h_2}$.
By \ref{getlevelling} applied to the blockade $(D_i:i\in J_1\cup \{h_2\})$, taking $H=H_3\cup \{h_2\}$, 
there exists $J_2\subseteq J_1\setminus H_3$ with $|J_2|\ge (|J_1|+1)/\rho-1$, and $E_j\subseteq D_j$ with $|E_j|\ge |D_j|/(4\rho)$ 
for all $j\in J_2$, and a $\left(D_j:j\in H_3\cup \{h_2\}\right)$-rainbow
levelling $\mathcal{M}=(M_0\LL M_m)$, with apex $v$, reaching $\mathcal{E}=(E_j:j\in J_2)$. (See figure \ref{fig:getbilevel2}.)
Thus $\mathcal{M}$ has height at most $\rho-1$.

This is almost what we want. There are two things to fix: there might be edges between $V(P)$ and $\bigcup_{j\in J_2}E_j$,
and there might be edges between  $V(P)\setminus \{u\}$ and $L_1\cupcup L_{t}$. To handle the first, we just use the bound on linkage, 
as follows. For each $j\in J_2$, and each $w\in V(P)$, at most $\lambda|A_j|$ vertices in $A_j$ are adjacent to $w$
(and none are adjacent to $u$.) Since $|V(P)\setminus \{u\}|\le \rho$, there are at most $\lambda\rho |A_j|$ vertices in $E_j$
with a neighbour in $V(P)$; and so there exists $F_j\subseteq E_j$ with $|F_j|\ge |E_j|-\lambda\rho |A_j|$, anticomplete to $V(P)$.
Since 
$$|E_j|\ge \frac{1}{4\rho}|D_j|\ge  \frac{9}{512\rho|I|}|A_j|,$$ 
it follows that 
$$|F_j|\ge \left(\frac{9}{512\rho|I|}-\lambda\rho\right)|A_j|\ge \frac{1}{64\rho|I|}|A_j|$$
since $\lambda\rho\le 1/(512\rho|I|)$.

Now we handle edges between  $V(P)\setminus \{u\}$ and $L_1\cupcup L_{t}$.
Let $j\in J_2$. Since $\mathcal{L}$ grades $\mathcal{C}$ forwards,
there exists $Y_j\subseteq L_{t}$ that covers 
$\bigcup_{i\in J_1,i\ge j}C_i$ and hence $\bigcup_{i\in J_1,i\ge j}F_i$, and 
is anticomplete to $\bigcup_{i\in J_1, i<j}C_i$. In particular, $Y_j$ is anticomplete to $\bigcup_{i\in H_3}C_i$.

Let $f\in F_j$. There is an induced path $Q$ between $f$ and $u$, with vertex set consisting of one 
vertex in each of $L_0, L_1\LL L_{t-1}, Y_j,\{f\}$. Since $u\in V(P)$, there is a subpath $Q'$ of $Q$ with one end $f$, and a subpath $P'$
of $P\setminus \{u\}$ with one end $v$, such that the subgraph induced on $V(P')\cup V(Q')$ is an induced path between $f,v$. Choose some such pair
$(P',Q')$ for each $f\in F_j$; we call $(|P'|,|Q'|)$ the {\em type} of $f$. Since $2\le |P'|<|P|\le \rho-1$, and 
$2\le |Q'|\le \rho+1$
(because $f$ has no neighbour in $V(P)$), there are
at most $\rho^2$ possible types. Hence there exists $G_j\subseteq F_j$ with $|G_j|\ge |F_j|/\rho^2$ such that all vertices in
$G_j$ have the same type. We call this common type the {\em type} of $j$. There exists $J_3\subseteq J_2$ with $|J_3|\ge |J_2|/\rho^2$
such that all $j\in J_3$ have the same type; let this common type be $(a,b)$. 
Let $P'$ be the subpath of $P$ with $a$ vertices and with one end $v$, and let $w$ be its other end. Let $N_0:=\{w\}$, and let
$N_1$ be the set of 
neighbours of $w$ in $L_{t-b+2}$. It follows that for each $j\in J_3$, and for each $f\in G_j$, there is an 
induced path $Q$ between $f$ and some vertex in $N_1$, with one   
vertex in each of $L_{t-b+2},L_{t-b+3}\LL L_{t-1}, Y_j$ and $\{f\}$; and the only edge between $V(P')$ and $V(Q')$ is the edge between $w$
and $N_1$. For $i=2\LL b-1$, let $N_i$ be the set of all vertices in $L_{t-b+i+1}$ that have no neighbour in $V(P')$ and have a neighbour in
$N_{i-1}$. Thus $\mathcal{N}$ is a levelling of height $b-1$ that  grades $\mathcal{G}=(G_j:j\in J_3)$.

Let the vertices of $P'$ in order be $w=p_1, p_2\LL p_{a}=v$. Then 
$$\mathcal{M}'=\left(\{p_1\},\{p_2\}\LL \{p_a\}=M_0,M_1\LL M_m\right)$$
is a levelling of height $a+m-1$, reaching $\bigcup_{j\in J_3}G_j$. Thus $(\mathcal{N},\mathcal{M}',\mathcal{G})$ is a bi-levelling.
Its height is $(a+m-1)+(b-1)\le 3\rho-3$, and so it 
satisfies the theorem.  This proves \ref{getbilevel}.~\bbox

\section{Selective covering}

Let $\mathcal{B}=(B_i:i\in I)$ be a blockade in $G$, and let $A\subseteq V(G)\setminus V(\mathcal{B})$ cover $V(\mathcal{B})$.
We wish to find a subset of $A$ that covers a significant amount and misses a linear fraction of several of the blocks of $\mathcal{B}$, assuming that $\mathcal{B}$ is sufficiently long.
That is the content of the next result.

\begin{thm}\label{selective}
Let $K\ge k\ge 1$ be integers, let $c>0$ with $K\ge (2+1/c)(k-1)$, let $0<\vare,\alpha\le 1$, let 
$\mathcal{B}=(B_i:i\in I)$ be a blockade 
of length $K$ in a graph $G$, and let $A\subseteq V(G)\setminus V(\mathcal{B})$ cover $V(\mathcal{B})$. Suppose that 
for each $i\in I$, the max-degree from $A$ to $B_i$ is less than 
$\vare|B_i|$. Suppose also that there is no
partition $\mathcal{P}$ of $A$ into at most $K^k$ sets, such that for each $X\in \mathcal{P}$ 
there exists $J\subseteq I$ with $|J|=k$ where $|N(X)\cap B_j|< K^k\alpha|B_j|$ for each $j\in J$. Then
there exists $X\subseteq A$ and a subset $J\subseteq I$ with $|J|=k$, such that 
$$|G|^{-c}\alpha\le |N(X)\cap B_i|/|B_i|< K^k\alpha+\vare$$ 
for each $i\in J$.
\end{thm}
\Proof
Let $\mathcal{J}$ be the set of all subsets of $I$ with cardinality $k$. For each $J\in \mathcal{J}$, choose $X(J)\subseteq A$ 
and $Y(J)\subseteq V(\mathcal{B})$ with the following properties:
\begin{itemize}
\item the sets $X(J)\;(J\in \mathcal{J})$ are pairwise disjoint subsets of $A$;
\item the sets $Y(J)\;(J\in \mathcal{J})$ are pairwise disjoint subsets of $V(\mathcal{B})$;
\item for each $J\in \mathcal{J}$, $Y(J)\subseteq N(X(J))\cap \bigcup_{j\in J}B_j$;
\item for each $J\in \mathcal{J}$ and each $j\in J$, $N(X(J))\cap B_j\subseteq \bigcup_{J'\in \mathcal{J}} Y(J')$;
\item for each $J\in \mathcal{J}$, and all distinct $i,j\in J$, $|Y(J)\cap B_i|/|B_i|\ge |G|^{-c}|Y(J)\cap B_j|/|B_j|$;
\item for each $J\in \mathcal{J}$, and each $j\in J$, $|Y(J)\cap B_j|< \alpha|B_j|$; and
\item subject to these conditions, $\bigcup_{J\in \mathcal{J}}X(J)$ is maximal.
\end{itemize}
This is possible, since we may set $X(J)=Y(J)=\emptyset$ to satisfy the first six conditions.
Suppose that $\bigcup_{J\in \mathcal{J}}X(J)=A$. Then the sets $X(J)\;(J\in \mathcal{J})$ form a partition of $A$ into at most $K^k$
subsets. Moreover, for each $J\in \mathcal{J}$, and each $j\in J$, we have $N(X(J))\cap B_j\subseteq \bigcup_{J'\in \mathcal{J}} Y(J')$,
and $|Y(J')\cap B_j|/|B_j|< \alpha$ for each $J'\in \mathcal{J}$, and so $|N(X(J))\cap B_j|< K^k\alpha|B_j|$, contrary to the hypothesis.
It follows that $\bigcup_{J\in \mathcal{J}}X(J)\ne A$.

Let $Y=\bigcup_{J\in \mathcal{J}}Y(J)$.
Choose $a\in A\setminus \bigcup_{J\in \mathcal{J}}X(J)$. For each $i\in I$, let $n_i$ be the number of neighbours of $a$
in $B_i\setminus Y$. 
Without loss of generality we may assume that $I=\{1\LL K\}$ where 
$$\frac{n_1}{|B_1|}\le \frac{n_2}{|B_2|}\le \cdots\le \frac{n_{K}}{|B_{K}|}.$$
If $n_k=0$, then we may add $a$ to $X(J)$, where $J=\{1\LL k\}$, contrary to the maximality of $\bigcup_{J\in \mathcal{J}}X(J)$.
Thus $n_k/|B_k|\ge |G|^{-1}$.
Choose an integer $t\ge 1$ with $t(k-1)+1\le K$, maximum such that 
$$\frac{n_{t(k-1)+1}}{|B_{t(k-1)+1}|}\ge |G|^{-1+(t-1)c}.$$

Since $n_{t(k-1)+1}/|B_{t(k-1)+1}|< \vare\le 1$, it follows that $|G|^{-1+(t-1)c}<1$, and so $-1+(t-1)c<0$, that is,
$t+1<2+1/c$. But $(2+1/c)(k-1)\le K$, and therefore 
$(t+1)(k-1)\le  K$. From the maximality of $t$, it follows that 
$$\frac{n_{(t+1)(k-1)+1}}{|B_{(t+1)(k-1)+1}|}< |G|^{-1+tc}.$$
Let $J:=\{j:t(k-1)+1\le j\le (t+1)(k-1)+1\}$, so $|J|=k$ and $J\in \mathcal{J}$.
Since 
$$\frac{n_{t(k-1)+1}}{|B_{t(k-1)+1}|}\ge |G|^{-1+(t-1)c},$$
it follows that 
$$\frac{n_{(t+1)(k-1)+1}}{|B_{(t+1)(k-1)+1}|}< |G|^{-1+tc}\le |G|^{c}\frac{n_{t(k-1)+1}}{|B_{t(k-1)+1}|}.$$
Consequently 
$n_{i}/|B_{i}|> |G|^{-c}n_{j}/|B_{j}|$ for all distinct $i,j\in J$.

For each $i\in I$, we define $N_i=N(a)\cap (B_i\setminus Y)$ (so $n_i=|N_i|$).
Define $X'(J)=X(J)\cup \{a\}$, and $Y'(J)=Y(J)\cup \bigcup_{j\in J}N_j$, and define $X'(J')=X(J')$ and $Y'(J')=Y(J')$
for all $J'\in \mathcal{J}\setminus \{J\}$.
From the maximality of $\bigcup_{J\in \mathcal{J}}X(J)$, replacing $X(J)$ by $X'(J)$ and $Y(J)$ by $Y'(J)$ 
violates one of the first six of the seven bullets above. The first four remain satisfied, so let us examine 
the fifth and sixth bullets.

Let $i,j\in J$ be distinct. Then $|Y(J)\cap B_i|/|B_i|\ge |G|^{-c}|Y(J)\cap B_j|/|B_j|$; and $|N_i|/|B_{i}|\ge  |G|^{-c}|N_j|/|B_{j}|$; and since
$|Y'(J)\cap B_i|=|Y(J)\cap B_i|+|N_i|$ and $|Y'(J)\cap B_j|=|Y(J)\cap B_j|+|N_j|$, it follows that 
$|Y'(J)\cap B_i|/|B_i|\ge |G|^{-c}|Y'(J)\cap B_j|/|B_j|$. Thus the fifth bullet remains satisfied. 

Consequently the sixth bullet is violated, and so there exists $j\in J$ such that $|Y'(J)\cap B_j|\ge \alpha|B_j|$. We claim that setting $X=X'(J)$ satisfies the theorem; and so we must check that 
$$|G|^{-c}\alpha\le |N(X'(J))\cap B_i|/|B_i|< K^k\alpha+\vare$$
for each $i\in J$. To prove these two inequalities, let $i\in J$.

Since $Y'(J)\subseteq N(X(J))\cap \bigcup_{j\in J}B_j$, and therefore $Y'(J)\cap B_i\subseteq N(X'(J))$, it follows that
$$|N(X'(J))\cap B_i|/|B_i|\ge |Y'(J)\cap B_i|/|B_i|\ge |G|^{-c}|Y'(J)\cap B_j|/|B_j|\ge  |G|^{-c}\alpha$$
(since the fifth bullet still holds).
This proves the first inequality. For the second, since
$|Y(J)\cap B_i|< \alpha|B_i|$ and 
$|N(a)\cap B_i|\le \vare|B_i|$, it follows that 
$|Y'(J)\cap B_i|< (\alpha+\vare)|B_i|$. But $|Y'(J')\cap B_i|< \alpha|B_i|$ for each $J'\in \mathcal{J}\setminus \{J\}$, and
$N(X'(J))\cap B_i\subseteq \bigcup_{J'\in \mathcal{J}} Y(J')$, and consequently $|N(X'(J))\cap B_i|\le (K^k\alpha+\vare)|B_i|$.
This proves the second inequality, and so 
proves \ref{selective}.~\bbox

We would like to obtain a better version of \ref{getbilevel}, where we can prescribe the height of 
the bi-levelling exactly. We show next that we can increase it by one, with the aid of 
\ref{selective}. 

\begin{thm}\label{extending}
Let $K\ge k\ge 1$ be integers, let $\gamma,\delta,\lambda, \eta, c>0$ with $K\ge (2+1/c)k$, and let $\mathcal{A}$ be a  
blockade with linkage at most $\lambda$ in a graph $G$, that is not $(\gamma,\delta)$-divergent,
where 
$\eta/2\ge K^{k+1}\delta|G|^c+\lambda+ \gamma$. Let
$(\mathcal{L}, \mathcal{M}, \mathcal{C})$ be an $\mathcal{A}$-rainbow bi-levelling with length $K$ and height $L$ say, 
where $\mathcal{C}$ has 
$\mathcal{A}$-size at least $\eta$.
Then there is an
$\mathcal{A}$-rainbow bi-levelling $(\mathcal{L}', \mathcal{M}', \mathcal{C}')$ with length $k$ and height $L+1$, 
such that $\mathcal{C}'$ is $\mathcal{C}$-rainbow, and has $\mathcal{C}$-size
at least $1/2$.
\end{thm}
\Proof
Let $\mathcal{L}=(L_0\LL L_\rho)$, and let $\mathcal{M}=(M_0\LL M_m)$. Thus $\rho+m=L$. 
Since $\mathcal{L}$ grades $\mathcal{C}$ forwards and $\mathcal{C}$ has length $K$, we may assume without loss of generality 
that $\mathcal{C}=(C_1\LL C_K)$, and for 
$1\le i\le K$ there exists $Y_i\subseteq L_{\rho}$ that covers $C_i\cupcup C_K$ and is anticomplete to $C_1\cupcup C_{i-1}$.

Now $M_m$ covers $V(\mathcal{C})$.
Moreover, for each $a\in M_m$, and each block $C$ of $\mathcal{C}$, since $\mathcal{A}$ has linkage at most $\lambda$, 
and $C$ has $\mathcal{A}$-size at
least $\eta$, it follows that $a$ has at most $(\lambda/\eta)|C|$ neighbours in $C$. 
\\
\\
(1) {\em If there is a partition $\mathcal{P}$ of $M_m$ into at most $K^{k+1}$ sets, such that for each $X\in \mathcal{P}$
there exists $J\subseteq \{1\LL K\}$ with $|J|=k+1$ such that $|N(X)\cap C_j|< K^{k+1}(\delta|G|^c/\eta)|C_j|$ for each $j\in J$, then there is a 
bi-levelling satisfying the theorem.}
\\
\\
Since $M_m$ covers $C_1$,
there exists $X\in \mathcal{P}$ that covers a subset $M\subseteq C_1$ with $\mathcal{C}$-size at least $K^{-k-1}$ and hence with
$\mathcal{A}$-size at least $\eta K^{-k-1}$. Since $X\in \mathcal{P}$, there exists $J'\subseteq \{1\LL K\}$ with $|J'|=k+1$ such that
$|N(X)\cap C_j|< K^{k+1}(\delta|G|^c/\eta)|C_j|$ for each $j\in J'$. Choose $J\subseteq J'$ with $1\notin J$, and with $|J|=k$.
For each $j\in J$ let $D_j:=C_j\setminus (N(X)\cap C_j)$; thus $D_j$ has $\mathcal{C}$-size at least 
$1-K^{k+1}(\delta|G|^c/\eta)$.
Since $M$ has $\mathcal{A}$-size at least $\eta K^{-k-1}\ge \delta$
and $\mathcal{A}$ is not $(\gamma,\delta)$-divergent, at most $\gamma|A_j|\le (\gamma/\eta)|C_j|$ vertices in $A_j$ have no neighbour in $M$. 
Let 
$C_j'$ be the set of vertices in $D_j$ that have a neighbour in $M$; thus $C_j'$ has $\mathcal{C}$-size at least 
$1-K^{k+1}(\delta|G|^c/\eta)- \gamma/\eta\ge 1/2$.
Let 
\begin{align*}
\mathcal{L}'&:=\left(L_0\LL L_{\rho-1}, \bigcup_{j\in J}Y_j\right)\\
\mathcal{M}'&:=\left(M_0\LL ,M_{m-1}, X, M\right)\\
\mathcal{C}'&:=\left(C_j':j\in J\right);
\end{align*}
then $(\mathcal{L}', \mathcal{M}', \mathcal{C}')$ is a 
bi-levelling satisfying the theorem. This proves (1).

\bigskip

By (1) we may therefore assume that there is no such partition $\mathcal{P}$.
By \ref{selective} (with $\vare$ 
replaced by $\lambda/\eta$, and $k$ replaced by $k+1$, and $\alpha$ replaced by $\delta|G|^c/\eta$) there exists $X\subseteq M_m$ and a 
subset $J'\subseteq \{1\LL K\}$ with $|J'|=k+1$, such that
$$\frac{\delta}{\eta}\le \frac{|N(X)\cap C_i|}{|C_i|}< K^{k+1}\frac{\delta|G|^c}{\eta}+\frac{\lambda}{\eta}$$ 
for each $i\in J'$.
Let $j_0$ be the smallest member of $J'$, and let $J:=J'\setminus \{j_0\}$.
Let 
$M:=N(X)\cap C_{j_0}$; thus $M$ has $\mathcal{C}$-size at least $\delta/\eta$, and so has
$\mathcal{A}$-size at least $\delta$.
For each $j\in J$, let $D_j$ be the set of vertices in $C_j$ that have no neighbour in $X$; thus $D_j$ has $\mathcal{C}$-size at least
$1-K^{k+1}(\delta|G|^c/\eta)-\lambda/\eta$.
Since the set of vertices in 
$C_j$ with no neighbour in $M$ has $\mathcal{A}$-size at most $\gamma$ and hence $\mathcal{C}$-size at most $\gamma/\eta$, 
it follows that $C_j'$ has $\mathcal{C}$-size at least
$1-K^{k+1}(\delta|G|^c/\eta)-\lambda/\eta- \gamma/\eta\ge 1/2$, where $C_j'$ is the set of vertices in $D_j$ with a neighbour in $M$.
Let
\begin{align*}
\mathcal{L}'&:=\left(L_0\LL L_{\rho-1}, \bigcup_{j\in J}Y_j\right)\\
\mathcal{M}'&:=\left(M_0\LL ,M_{m-1}, X, M\right)\\
\mathcal{C}'&:=\left(C_j':j\in J\right);
\end{align*}
then $(\mathcal{L}', \mathcal{M}', \mathcal{C}')$ is a
bi-levelling satisfying the theorem. This proves \ref{extending}.~\bbox

By combining \ref{getbilevel} and \ref{extending}, we obtain a version of \ref{getbilevel} where we can specify the height of 
the bi-levelling exactly, the following.

\begin{thm}\label{exactbilevel}
Let $k\ge 1$ be an integer, and let $c>0$. Let $\rho:=\lceil 1+1/c\rceil$. Let $\ell\ge 3\rho-2$ be an integer.
Let $K:= \lceil k(3+1/c)^{\ell+2}\rceil$. Let $\gamma,\lambda>0$
with $\lambda,\gamma\le 2^{-8-\ell}/(\rho^3K)$.
Let $G$ be a graph and define $\delta:=K^{-K}|G|^{-c}$.
Let $\mathcal{A}=(A_i:i\in I)$ be 
a blockade in $G$ of length $K$, with linkage at most $\lambda$, that is not  $(\gamma,\delta)$-divergent.
Then there is an $\mathcal{A}$-rainbow bi-levelling $(\mathcal{L}, \mathcal{M}, \mathcal{C})$ with length $k$ and height
$\ell$, such that $\mathcal{C}$ has $\mathcal{A}$-size at least $2^{4-\ell}/(\rho^3K)$.
\end{thm}
\Proof
We would like to apply \ref{getbilevel} with $k$ replaced by $\lfloor k(3+1/c)^{\ell-2}\rfloor$. 
We must check that 
\begin{eqnarray*}
\lfloor k(3+1/c)^{\ell-2}\rfloor&\le& K/\rho^4\\
\lambda&\le& 1/(512\rho^2K)\\
\gamma&\le& 3/(256K) \\
\delta&\le& 3\rho/(128 K^2)\\
(256K\delta/3)^{\rho-1}|G| &\le& 1.
\end{eqnarray*}
The first follows since since 
$$k(3+1/c)^{\ell-2}\le  k(3+1/c)^{\ell+2}/\rho^4\le K/\rho^4$$
(because  $\rho\le 2+1/c$).
The second and third are implied by the hypothesis $\lambda,\gamma\le 2^{-8-\ell}/(\rho^3K)$, since $\rho\ge 2$ and therefore 
$\ell\ge 4$ (because 
$\ell\ge 3\rho-2$). The fourth follows since  
$$\delta=K^{-K}|G|^{-c}\le K^{-K}\le 3/(128K^2)\le 3\rho/(128 K^2)$$
(because $K\ge 3^4$, and so $K^{2-K}\le 3/128$). The fifth follows since
$$(256K\delta/3)^{\rho-1}\le (\delta K^K)^{\rho-1}=|G|^{-c(\rho-1)}\le |G|^{-1}$$
(because $256K/3\le K^K$). 
Thus we can apply \ref{getbilevel}. We deduce that
there is an $\mathcal{A}$-rainbow bi-levelling $(\mathcal{L}', \mathcal{M}', \mathcal{C}')$ with length $\lfloor k(3+1/c)^{\ell-2}\rfloor$ and height
at most $3\rho-3$, such that $\mathcal{C}$ has $\mathcal{A}$-size at least $1/(64\rho^3K)$.

Let its height be $\ell-t$; thus $1\le t\le \ell-2$. Define $K_{0}=k$, and for $i=1\LL t$ let $K_{i}:=\lceil(2+1/c)K_{i-1}\rceil$.
Thus $K_i\le \lfloor (3+1/c)K_{i-1}\rfloor$, and so 
$$K_t\le \lfloor k(3+1/c)^{t}\rfloor \le \lfloor k(3+1/c)^{\ell-2}\rfloor\le  K/\rho^4.$$
For $i=0\LL t$ define $\eta_i=2^{i-t-6}/(\rho^3K)$.
Thus (setting $\mathcal{L}_{t}=\mathcal{L}'$ and $\mathcal{M}_{t}=\mathcal{M}'$, and letting $\mathcal{C}_{t}$ be a sub-blockade
of $\mathcal{C}'$ with the right length) we deduce that there is an
$\mathcal{A}$-rainbow bi-levelling  $(\mathcal{L}_{t}, \mathcal{M}_{t}, \mathcal{C}_{t})$ with length $K_{t}$ and height $\ell-t$,
such that $\mathcal{C}_{t}$ has $\mathcal{A}$-size
at least $\eta_{t}$.
\\
\\
(1) {\em $\eta_s/8\ge K_s^{K_{s-1}+1}\delta|G|^c+\lambda+ \gamma$ for $1\le s\le t$.}
\\
\\
Certainly $\eta_s = 2^{s-t-6}/(\rho^3K)\ge 2^{-t-5}/(\rho^3K)$. 
Moreover,   
$2^{\ell+7}\le 3^{2\ell+1}$ since $\ell\ge 4$; so 
$$2^{\ell+7}(2+1/c)^{3}\le (3+1/c)^{2\ell+4}\le  K^{2}.$$
Consequently
$$\frac{\eta_s}{16}\ge \frac{1}{2^{t+9}\rho^3K}\ge \frac{1}{2^{\ell+7}\rho^3K}\ge \frac{1}{2^{\ell+7}(2+1/c)^{3}K}\ge \frac{1}{K^3}.$$
Since $K-3\ge K_{s-1}+1$, and therefore $K^{K-3}\ge K^{K_{s-1}+1}\ge K_s^{K_{s-1}+1}$, it follows that
$$\eta_s/4\ge K^{-3} \ge K_s^{K_{s-1}+1}K^{-K}=K_s^{K_{s-1}+1}\delta|G|^c.$$
But also, since $\lambda,\gamma\le 2^{-8-\ell}/(\rho^3K)$, it follows that
$$\frac{\eta_s}{16}\ge 2^{-t-9}/(\rho^3K)\ge 2^{-7-\ell}/(\rho^3K)\ge \gamma+\lambda.$$
Adding, we deduce that 
$$\eta_s/8\ge  K_s^{K_{s-1}-1}\delta|G|^c+ \gamma+\lambda.$$
This proves (1).

\bigskip

By (1), taking $s=t$, we may apply \ref{extending}, replacing $\eta, K, k$ by $\eta_{t}, K_t$, and  
$K_{t-1}$ respectively. We deduce that there is an
$\mathcal{A}$-rainbow bi-levelling $(\mathcal{L}_{t-1}, \mathcal{M}_{t-1}, \mathcal{C}_{t-1})$ with length $K_{t-1}$ and 
height $\ell-(t-1)$,
such that $\mathcal{C}_{t-1}$ is $\mathcal{A}$-rainbow and $\mathcal{C}_{t-1}$ has $\mathcal{A}$-size
at least $\eta_{t-1}$.

Choose $s\le t-1$ with $s\ge 0$ minimum such that there is an
$\mathcal{A}$-rainbow bi-levelling  $(\mathcal{L}_{s}, \mathcal{M}_{s}, \mathcal{C}_{s})$ with length $K_{s}$ and height $\ell-s$,
such that $\mathcal{C}_{s}$ has $\mathcal{A}$-size
at least $\eta_s$.
Suppose that $s>0$. By (1) we may apply \ref{extending},
replacing $\eta$ by $\eta_{s}$, and replacing $K$ by $K_s$ and replacing $k$ by $K_{s-1}$, giving a contradiction to the minimality of $s$. 

Thus $s=0$.
Hence there is an
$\mathcal{A}$-rainbow bi-levelling  $(\mathcal{L}_{0}, \mathcal{M}_{0}, \mathcal{C}_{0})$ with length $K_{0}=k$ and height $\ell$,
such that $\mathcal{C}_{0}$ has $\mathcal{A}$-size
at least $\eta_0\ge 2^{4-\ell}/(\rho^3K)$.
This proves \ref{exactbilevel}.~\bbox

\section{A digression}
The result \ref{exactbilevel} needs to be strengthened further for its use in this paper, but as it stands it is 
already quite strong. For instance, it gives an improvement over a result of \cite{pure5} which was one of the 
main theorems of that paper.
In \cite{pure5} we proved:

\begin{thm}\label{pure5thm}
Let $c>0$ with $1/c$ an integer, and let $\ell\ge 4/c+5$ be an integer.
Then there
exists $\vare>0$ such that
if $G$ is an $\vare$-sparse $(\vare|G|^{1-c}, \vare|G|)$-coherent graph with $|G|>1$, then $G$ has an induced cycle of length $\ell$.
\end{thm}

With the aid of \ref{exactbilevel} we can do a little better:
\begin{thm}\label{betterpure5thm}
Let $c>0$ with $1/c$ an integer, and let $\ell\ge 3/c+3$ be an integer.
Then there
exists $\vare>0$ such that
if $G$ is an $\vare$-sparse $(\vare|G|^{1-c}, \vare|G|)$-coherent graph with $|G|>1$, then $G$ has an induced cycle of length $\ell$.
\end{thm}
\Proof
Let $K:= (3+1/c)^{\ell}$, and $\rho=1+1/c$.
Let $\vare>0$ satisfy $2K\vare\le K^{-K}$ and $2K\vare\le 2^{-6-\ell}/(\rho^3K)$.
We claim that $\vare$ satisfies the theorem.

Let $G$ be  an $\vare$-sparse $(\vare|G|^{1-c}, \vare|G|)$-coherent graph.
By \ref{big}, $|G|>1/\vare\ge K$, and so $\lfloor |G|/K\rfloor \ge |G|/(2K)$; and consequently
there is a blockade $\mathcal{A}$ in $G$ of length $K$ and width at least $|G|/(2K)$. Its linkage is at most
$2K\vare$ since $G$ is $\vare$-sparse, and it is not $(2K\vare|G|^{-c}, 2K\vare)$-divergent since $G$ is  
$(\vare|G|^{1-c}, \vare|G|)$-coherent. 

Let $\vare':=\gamma=2K\vare$; so $\mathcal{A}$ has linkage at most $\vare'$.
Let $\delta:=K^{-K}|G|^{-c}$; then $\delta\ge 2K\vare|G|^{-c}$ since $K^{-K}\ge 2K\vare$,
and so $\mathcal{A}$ is not $(\gamma,\delta)$-divergent.
Let $\ell':=\ell-2$.
Thus $\ell'\ge 3\rho-2$ since $\ell\ge 3/c+3$.
Also, $\vare',\gamma\le 2^{-8-\ell'}/(\rho^3K)$.
By \ref{exactbilevel} with $k=1$, there is an $\mathcal{A}$-rainbow bi-levelling $(\mathcal{L}, \mathcal{M}, \mathcal{C})$ 
with length $1$ and height
$\ell'$. Choose $w$ in the unique block of $\mathcal{C}$. Then $w$ has a neighbour $u$ in the base of $\mathcal{L}$ and a 
neighbour $v$
in the base
of $\mathcal{M}$, and there is an induced path of length $\ell'$ between $u,v$ whose internal vertices are anticomplete to 
the block of $\mathcal{C}$; and so adding $w$ to this path gives an induced cycle of length $\ell$. This proves \ref{betterpure5thm}.~\bbox

\section{Bi-gradings}

Let $(\mathcal{L}, \mathcal{M},\mathcal{C})$ be a bi-levelling. Thus $\mathcal{L}$ grades $\mathcal{C}$ forwards, but 
$\mathcal{M}$ does not, and next we want to arrange that $\mathcal{M}$ also grades $\mathcal{C}$. We can do this with the same 
argument
that we used for $\mathcal{L}$, and that gives a corresponding ordering of the boxes of $\mathcal{C}$ (or rather, of the contraction of a sub-blockade of 
$\mathcal{C}$ that survives this argument), but the two orderings might be very similar, and that turns out not to be useful.
What we need is that $\mathcal{M}$ grades $\mathcal{C}$ in the opposite order from $\mathcal{L}$, and that is the subject of this section.

If we start
with a blockade with sufficient length,
the result \ref{exactbilevel} provides us with a bi-levelling of any desired height, and any desired length,
just at the cost of shrinking the blocks by constant factors.
But to persuade the part $\mathcal{M}$ of the bi-levelling to grade $\mathcal{C}$ backwards, we no longer have the luxury of 
linear shrinking; now we will have to 
shrink the blocks by fractions that are polynomial in $|G|$. (This is why the proof only proves \ref{mainthm} and not \ref{conj}.)

Let us say this more precisely. Let $\mathcal{M}=(M_0\LL M_m)$ be a levelling and $\mathcal{B}$ a blockade.
We say that $\mathcal{M}$ {\em grades $\mathcal{B}$ backwards} if $\mathcal{M}$ 
reaches $\bigcup_{i\in I}B_i$, and
for each $j\in I$ there exists
$Y\subseteq M_m$ that
covers $\bigcup_{i\in I, i\le j}B_i$ and
is anticomplete to $\bigcup_{i\in I, i>j}B_i$.
We say that a bi-levelling $(\mathcal{L}, \mathcal{M},\mathcal{C})$ is a {\em bi-grading} if $\mathcal{M}$ grades
$\mathcal{C}$ backwards. Other definitions (length, height, $\mathcal{A}$-rainbow) are the same as for a bi-levelling.

\begin{thm}\label{bigrade}
Let $\ell,k\ge 1$ be integers, and let $0<c,d,\lambda'<1$, such that $\ell\ge 3\lceil 1/c\rceil+1$.
Then there exist an integer $K\ge 1$ and $\lambda>0$ with the following property.
Let $G$ be a graph, and let $\mathcal{A}=(A_i:i\in I)$ be
a blockade in $G$ of length $K$, with linkage at most $\lambda$, that is not  $(2^{-6-2\ell}/K,K^{-K}|G|^{-c})$-divergent.
Then there is an $\mathcal{A}$-rainbow bi-grading $(\mathcal{L}, \mathcal{M}, \mathcal{B})$ with length $k$ and height
$\ell$, such that $\mathcal{B}$ has $\mathcal{A}$-size at least
$2^{-k-2\ell}K^{-1-K}|G|^{-d}$
and linkage at most $\lambda'$.
\end{thm}




\Proof Let $K_k:=1$, and for $t=k-1,k-2\LL 0$ let $K_t:= \lceil(2+1/d)K_{t+1}+1\rceil$.
Let $K:= \lceil K_0(3+1/c)^{\ell+2}\rceil$, and let
$\rho=\lceil1+1/c\rceil$; then $\ell\ge 3\rho-2$, since $\ell\ge 3\lceil 1/c\rceil+1$.
It follows that $\rho^3\le (\ell+2)^3/27\le 2^{\ell-2}$, since $\ell\ge 7$.
Let $\eta:=2^{6-2\ell}/K$, and $\lambda= \lambda'\eta/(k2^{k+11})$.
We claim that $K,\lambda$ satisfy the theorem.
Let $G$ be a graph and let $\mathcal{A}$ be
a blockade in $G$ of length $K$, with linkage at most $\lambda$, that is not
$(2^{-6-2\ell}/K,K^{-K}|G|^{-c})$-divergent.
Since $2^{-6-2\ell}/K\le  2^{-8-\ell}/(\rho^3K)$ and $\eta\le 2^{4-\ell}/(\rho^3K)$,
and
$$\lambda\le 2^{6-2\ell}/(Kk2^{k+11})\le 2^{-2\ell-6}/K \le 2^{-8-\ell}/(\rho^3K),$$
\ref{exactbilevel} implies that there is an $\mathcal{A}$-rainbow bi-levelling $(\mathcal{L}, \mathcal{M}, \mathcal{C}_0)$ with length $K_0$
and height
$\ell$, such that $\mathcal{C}_0$ has $\mathcal{A}$-size at least $\eta$, and hence has linkage
at most $\lambda/\eta$.
The theorem is invariant under re-indexing $\mathcal{A}$;
and so we may assume that $\mathcal{C}_0$ is a contraction of a sub-blockade of $\mathcal{A}$.
Let $\mathcal{A}=(A_i:i\in I)$, and let $\mathcal{C}_0=(C^0_i:i\in I_0)$, where $I_0\subseteq I$ and $C^0_i\subseteq A_i$ for each $i\in I_0$. Let $\mathcal{M}$ have height $m$, and let its base be $M_m$.

Suppose inductively that we have defined $i_1\LL i_{t}\in I_0$, and $I_0,I_1\LL I_{t}$, and $D_{i_1}\LL D_{i_t}$, and
$C^j_i$ for $0\le j\le t$ and each $i\in I_j$, with the following properties for $1\le j\le t$:
\begin{itemize}
\item $I_j\subseteq I_{j-1}$ with $|I_j|= K_j$, and $i_j\in I_{j-1}\setminus I_j$,  and $i>i_j$ for all $i\in I_j$ (and consequently $i_1<i_2<\cdots < i_{t}$);
\item $C^j_i\subseteq C^{j-1}_i$ and $|C^j_i|\ge |C^{j-1}_i|/2$ for all $i\in I_j$;
\item $D_{i_j}\subseteq C^{j-1}_{i_j}$ and has $\mathcal{A}$-size at least $2^{1-k-2\ell}K^{-1-K}|G|^{-d}$;
\item there exists $X\subseteq M_m$ such that $X$ covers $D_{i_j}$, and is anticomplete to $C^j_i$
for all $i\in I_j$;
\item for all $h\in \{1\LL j-1\}$ the max-degree from $D_{i_j}$ to $D_{i_h}$ is at most $\lambda'|D_{i_h}|/(4k)$;
\item for all $i\in I_j$, the max-degree from $C^j_i$ to $D_{i_j}$ is at most $\lambda'|D_{i_j}|/(4k)$.
\end{itemize}
If $t=k$ the inductive definition is complete, so we assume that $0\le t<k$; and now we need to choose $i_{t+1}$, $I_{t+1}$,
$D_{i_{t+1}}$,
and $C^{t+1}_i$ for each $i\in I_{t+1}$, so that the bullets are satisfied with $t$ replaced by $t+1$.
\\
\\
(1) {\em There exist $i_{t+1}\in I_t$, and a subset $I_{t+1}\subseteq I_t$ with $|I_{t+1}|=K_{t+1}$, and a subset $X\subseteq M_m$,
and a subset $D_{i_{t+1}}\subseteq C_{i_{t+1}}$ with $\mathcal{A}$-size at least $2^{1-k-2\ell}K^{-K}|G|^{-d}$,
such that $i>i_{t+1}$ for each $i\in I_{t+1}$,
and $X$ covers $ D_{i_{t+1}}$ and is anticomplete to at least $3/4$ of $C^t_i$ for each $i\in I_{t+1}$.}
\\
\\
For $i\in I_t$, the max-degree from $M_m$ to $A_i$ is at most
$$\lambda|A_i|\le \lambda |C^0_i|/\eta\le \lambda 2^{t}|C^t_i|/\eta.$$
Since $|I_t|=K_t\ge (2+1/d)K_{t+1}$, we can apply  \ref{selective} to $M_m$ and the blockade $(C^t_i:i\in I_t)$,
with $\vare, c, K, k,\alpha$ replaced by $\lambda 2^t/\eta, d, K_t, K_{t+1}+1,K^{-K}/8$ respectively. We deduce that either:
\begin{itemize}
\item there is a
partition $\mathcal{P}$ of $M_m$ into at most $K^K$ sets, such that for each $X\in \mathcal{P}$
there exists $J\subseteq I_t$ with $|J|=K_{t+1}+1$ such that $|N(X)\cap C^t_j|< |C^t_j|/8$ for each $j\in J$; or
\item there exists $X\subseteq M_m$ and a subset $J\subseteq I_t$ with $|J|=K_{t+1}+1$, such that
$$|G|^{-d}K^{-K}/8\le |N(X)\cap C_i^t|/|C_i^t|< 1/8+\lambda2^t/\eta$$
for each $i\in J$.
\end{itemize}

Suppose the first happens, and let $\mathcal{P}$ be such a partition. Let $i_{t+1}$ be the smallest member of $I_t$.
Since $M_m$ covers $C^t_{i_{t+1}}$, there exists $X\in \mathcal{P}$ such that
$$|N(X)\cap C^t_{i_{t+1}}|\ge K^{-K}|C^t_{i_{t+1}}|.$$
Define $D_{i_{t+1}}=N(X)\cap C^t_{i_{t+1}}$; thus $D_{i_{t+1}}$ has $\mathcal{C}_0$-size at least $K^{-K}2^{-t}\ge K^{-K}2^{1-k}$,
and so has $\mathcal{A}$-size at least $2^{1-k-2\ell}K^{-1-K}|G|^{-d}$, since $2^{1-k-2\ell}K^{-1-K}|G|^{-d}\le \eta K^{-K}2^{1-k}$.
There exists $J\subseteq I_t$ with $|J|=K_{t+1}+1$ such that $|N(X)\cap C^t_j|< |C^t_j|/8$ for each $j\in J$;
choose $I_{t+1}\subseteq J\setminus \{i_{t+1}\}$ with cardinality $K_{t+1}$. Since
$|N(X)\cap C^t_i|< |C^t_i|/8\le |C^t_i|/4$ for each $i\in I_{t+1}$, in this case (1) holds.

So we may assume the second happens, and so there exists $X\subseteq M_m$ and a subset $J\subseteq I_t$ with $|J|=K_{t+1}+1$, such that
$$|G|^{-d}K^{-K}/8\le |N(X)\cap C_i^t|/|C_i^t|< 1/8+\lambda2^t/\eta$$
for each $i\in J$. Let $i_{t+1}$ be the smallest element of $J$, and define $I_{t+1}=J\setminus \{i_{t+1}\}$. Define
$D_{i_{t+1}}=N(X)\cap C^t_{i_{t+1}}$; thus $D_{i_{t+1}}$ has $\mathcal{C}_0$-size at least  
$|G|^{-d}K^{-K}2^{-t-3}\ge |G|^{-d}K^{-K}2^{-k-2}$, and hence
has $\mathcal{A}$-size at least $2^{1-k-2\ell}K^{-1-K}|G|^{-d}$, since $2^{1-k-2\ell}K^{-1-K}|G|^{-d}\le \eta|G|^{-d}K^{-K}2^{-k-2}$.
Since $1/8+\lambda2^t/\eta\le 1/4$
it follows that again (1) holds. This proves (1).

\bigskip
Choose $i_{t+1}$, $I_{t+1}$, $X$ and $D_{i_{t+1}}$ as in (1). We claim this satisfies the conditions for the inductive definition.
For each $i\in I_{t+1}$, since the max-degree from $D_{i_{t+1}}$
to $C_i^t$ is at most $2^t(\lambda/\eta)|C_i^t|$, there are at most $|C_i^t|/4$ vertices in $C_i^t$ that have at least $2^{t+2}(\lambda/\eta)|D_{i_{t+1}}|$
neighbours in $D_{i_{t+1}}$. Consequently, at least half the vertices in  $C_i^t$ have fewer than $2^{t+2}(\lambda/\eta)|D_{i_{t+1}}|$
neighbours in $D_{i_{t+1}}$ and have no neighbour in $X$. Define $C_i^{t+1}$ to be the set of all such vertices.
For $1\le h\le t$, the max-degree from $D_{i_{t+1}}$ to $D_{i_h}$ is at most the max-degree from $C_{i_{t+1}}^t$ to $D_{i_h}$,
and hence is at most $\lambda'|D_{i_h}|/(4k)$. Thus, this completes the inductive definition.

For $1\le h\le k$, let $B_{i_h}$ be the set of all vertices in $D_{i_h}$ that have at most $(\lambda'/2)|D_{i_j}|$
neighbours in $D_{i_j}$ for all $j$ with $h<j\le k$. Since every vertex in $D_{i_j}$ has at most $(\lambda'/(4k))|D_{i_h}|$
neighbours in $D_{i_h}$, there are at most $|D_{i_h}|/(2k)$ vertices in $D_{i_h}$ that have at least $(\lambda'/2)|D_{i_j}|$ neighbours
in $D_{i_j}$; and so $|B_{i_h}|\ge |D_{i_h}|/2$. Thus for all $h,j\in \{1\LL t\}$ with $h<j$, the max-degree from $B_{i_h}$
to $B_{i_j}$ is at most $(\lambda'/2)|D_{i_j}|\le \lambda'|B_{i_j}|$; and the max-degree from $B_{i_j}$ to $B_{i_h}$ is at most
$(\lambda'/4k)|D_{i_h}|\le \lambda'|B_{i_h}|$. This proves \ref{bigrade}.~\bbox

We recall that the {\em shrinkage} of a blockade $\mathcal{B}=(B_i:i\in I)$ is the number $\sigma$ such that $|G|^{1-\sigma}$ 
is the width of $\mathcal{B}$. Let us recast \ref{bigrade} in terms of shrinkage.

\begin{thm}\label{bigradeshrink}
Let $\ell$ be an integer, and let $c, \sigma, \sigma'>0$ with $\sigma<\sigma'$ and $1\ge c> 1/\lfloor(\ell-1)/3\rfloor$.
Let $K'>0$ be an integer, and let $\lambda'>0$. Then 
there exist $\lambda>0$ and integers $N,K>0$ 
with the following property. Let $G$ be a graph with $|G|\ge N$, and let $\mathcal{A}$ be a blockade in $G$ of length $K$, with linkage at most $\lambda$ and shrinkage
at most $\sigma$, such that $\mathcal{A}$ is not $(|G|^{-c}, |G|^{-c})$-divergent. Then there is 
an $\mathcal{A}$-rainbow bi-grading $(\mathcal{L}, \mathcal{M}, \mathcal{B})$ with height
$\ell$, such that $\mathcal{B}$ has length $K'$,
linkage at most $\lambda'$ and shrinkage at most $\sigma'$.
\end{thm}
\Proof Choose $c'$ with $c> c'> 1/\lfloor(\ell-1)/3\rfloor$. Choose $d$ with $0<d<\sigma'-\sigma$. 
Choose $K, \lambda$ such that \ref{bigrade} is satisfied with $k$ replaced by $K'$. 
Choose $N\ge 0$ such that $N^c\ge 2^{6+2\ell}K$ and $N^{c-c'}\ge K^K$ and $N^{\sigma'-\sigma-d}\ge 2^{K'+2\ell}K^{1+K}$.
We claim that $N,K,\lambda$ satisfy the theorem.

Let $G$ be a graph with $|G|\ge N$, and let $\mathcal{A}$ be a blockade in $G$ of length $K$, with linkage at most $\lambda$ and shrinkage
at most $\sigma$, such that $\mathcal{A}$ is not $(|G|^{-c}, |G|^{-c})$-divergent.  Since $2^{-6-2\ell}/K\ge |G|^{-c}$ and
$K^{-K}|G|^{-c'}\ge |G|^{-c}$, it follows that $\mathcal{A}$ is not $(2^{-6-2\ell}/K,K^{-K}|G|^{-c'})$-divergent. 

By \ref{bigrade},
there is an $\mathcal{A}$-rainbow bi-grading $(\mathcal{L}, \mathcal{M}, \mathcal{B})$ with length $K'$ and height
$\ell$, such that $\mathcal{B}$ has $\mathcal{A}$-size at least
$2^{-K'-2\ell}K^{-1-K}|G|^{-d}$
and linkage at most $\lambda'$. Hence $\mathcal{B}$ has shrinkage at most $\sigma'$, since 
$$|G|^{1-\sigma'}\le |G|^{1-\sigma}2^{-K'-2\ell}K^{-1-K}|G|^{-d}$$
(because $N^{\sigma'-\sigma-d}\ge 2^{K'+2\ell}K^{1+K}$).
This proves \ref{bigradeshrink}.~\bbox

\section{Enforcement}

Let $\mathcal{B}=(B_i:i\in I)$ be a blockade in $G$, and let $H$ be a graph such that there is a $\mathcal{B}$-rainbow
copy $J$ of $H$. Since each vertex of $J$ belongs to some block of $\mathcal{B}$, all different, and the blocks are numbered 
by integers,
this defines an order on the vertices of $J$, and hence, via the isomorphism, an order on the vertices of $H$. We cared about this order in       
\cite{pure6}, but here the order does not concern us. What does concern us are the first and last vertices of $J$, and correspondingly of $H$.

Let $J$ be $\mathcal{B}$-rainbow, let $u\in V(J)$, and let $i\in I$ such that $u\in B_i$; we say that $u\in V(J)$ is  
the {\em $\mathcal{B}$-first} vertex of $J$ if there is no $h\in I$ with $h<i$ such that $B_h\cap V(J)\ne \emptyset$.
We define the {\em $\mathcal{B}$-last} vertex of $J$ similarly.

Now let $H$ be a graph, let $v\in V(H)$ and $\mathcal{B}$ be a blockade in $G$.
We say that a $\mathcal{B}$-rainbow copy $J$ of $H$ is {\em $v$-first} if there is an isomorphism
$\phi$ from $H$ to $J$, such that $\phi(v)$ is the $\mathcal{B}$-first vertex of $J$. We define {\em $v$-last} similarly.

Now let $H$ be a graph, let $K,N>0$ be integers, and let $0< \lambda,\sigma,c\le 1$. We say that
$(N,K,\sigma,\lambda,c)$ {\em forces} $H$ if for every $(|G|^{1-c},|G|^{1-c})$-coherent graph $G$ with $|G|\ge N$, 
and for every blockade $\mathcal{A}$ in $G$ of length $K$, 
shrinkage at most $\sigma$, and linkage at most $\lambda$, there is an $\mathcal{A}$-rainbow copy of $H$.
Similarly if $v\in V(H)$, we say that $(N,K,\sigma,\lambda,c)$ {\em forces $H$ $v$-first} if 
the same statement holds where the $\mathcal{A}$-rainbow copy of $H$ is $v$-first.
We define 
{\em forces $H$ $v$-last}, and {\em forces $H$ $u$-first and $v$-last} similarly.

From \ref{bigradeshrink} we deduce:
\begin{thm}\label{addhandle}
Let $H$ be a graph, obtained from a graph $H'$ by adding a handle of length $\ell$ with ends $u,v$. 
Let $N',K'\ge 1$ be integers, and let $0< c,\sigma,\sigma',\lambda'\le 1$, such that $\sigma<\sigma'$ and $c-\sigma> 1/\lfloor(\ell-1)/3\rfloor$,
and $(N',K',\sigma',\lambda',c)$ forces $H'$ $u$-first and $v$-last.
Then exist $\lambda>0$ and integers $N,K>0$ such that 
$(N,K,\sigma,\lambda,c)$ forces $H$.
\end{thm}
\Proof Since $c-\sigma> 1/\lfloor(\ell-1)/3\rfloor$, we can apply \ref{bigradeshrink} with $c$ replaced by $c-\sigma$; let 
$N,K,\lambda$ satisfy \ref{bigradeshrink} with $c$ replaced by $c-\sigma$. We may choose $N\ge N'$.
We claim that $N,K,\lambda$ satisfy the theorem.

Let $G$ be a $(|G|^{1-c},|G|^{1-c})$-coherent graph with $|G|\ge N$,
and let $\mathcal{A}=(A_i:i\in I)$
be a blockade in $G$ of length $K$, shrinkage at most $\sigma$ and linkage at most $\lambda$.
We must show that there is an $\mathcal{A}$-rainbow copy of $H$.
Since $\mathcal{A}$ has shrinkage at most $\sigma$, and  $G$ is $(|G|^{1-c},|G|^{1-c})$-coherent, 
it follows that $\mathcal{A}$ is not $(|G|^{\sigma-c}, |G|^{\sigma-c})$-divergent.
Since $N,K,\lambda$ satisfy \ref{bigradeshrink} with $c$ replaced by $c-\sigma$, there is
an $\mathcal{A}$-rainbow bi-grading $(\mathcal{L}, \mathcal{M}, \mathcal{B})$ with height
$\ell$, such that $\mathcal{B}$ has length $K'$,
linkage at most $\lambda'$ and shrinkage at most $\sigma'$.
Let $\mathcal{B}=(B_i:i\in I')$. 
Since  $(N',K',\sigma',\lambda',c)$ forces $H'$ $u$-first and $v$-last in $G$, there is a $u$-first $v$-last $\mathcal{B}$-rainbow copy 
$J$ of $H'$. Let $\phi$ be the corresponding isomorphism from $H$ to $J$, and let $\phi(u)\in B_h$ and $\phi(v)\in B_k$, 
where $h,k\in I'$. Thus for every $w\in V(H)\setminus \{u,v\}$, there exists $i\in I$ with $\phi(w)\in B_i$ and $h<i<k$.

Since $\mathcal{L}$ grades $\mathcal{B}$ forwards, there 
exists a subset $Y$ of the base of $\mathcal{L}$ that covers $B_k$ and is anticomplete to $V(J)\setminus B_k$;
and in particular,
there exists a vertex $y$ in the base of $\mathcal{L}$ that is adjacent to $\phi(v)$ and nonadjacent to all other vertices of 
$J$. 
Since $\mathcal{M}$
grades $\mathcal{B}$ backwards, similarly there is a vertex $x$ in the base of $\mathcal{M}$ that is adjacent to $\phi(u)$ and nonadjacent to all other vertices of $J$.
By \ref{bilevelpath}, there is an induced path $P$ between $x,y$, of length $\ell$, with interior
disjoint from and anticomplete to $V(J)\setminus \{\phi(u),\phi(v)\}$, and adding $P$ to $J$ gives an $A$-rainbow copy of $H$.
This proves \ref{addhandle}.~\bbox

\section{Covering by leaves}

We need to apply a theorem of \cite{pure6}, the following:

\begin{thm}\label{fullleafcover}
Let $K'\ge 0$ be an integer, and let $0<c, \sigma, \sigma',\lambda'\le 1$ with $\sigma<\sigma'<c$.
Then there exist $\lambda>0$ and integers $K,N>0$ with the following property.
Let $G$ be a $(|G|^{1-c},|G|^{1-c})$-coherent graph with $|G|\ge N$, and let $\mathcal{A}=(A_i:i\in I)$ be a
blockade of length $K$ in $G$, with shrinkage at most $\sigma$ and linkage at most $\lambda$.
Then there exists $I'\subseteq I$ with $|I'|=K'$, such that for every partition $(H,J)$ of $I'$,
there exists $B_h\subseteq A_h$ for each $h\in H$,
where
\begin{itemize}
\item $(B_h:h\in H)$ has shrinkage at most $\sigma'$ and linkage at most $\lambda'$; and
\item for all $h\in H$ and
all $j\in J$ there exists $X\subseteq A_j$ that covers $B_h$ and is anticomplete to $\bigcup_{i\in H\setminus \{h\}}B_{i}$.
\end{itemize}
\end{thm}

\begin{thm}\label{twoleaves}
Let $H$ be a graph, and let $u,v$ be distinct nonadjacent vertices of $H$, both with degree one in $H$.
Let $H' := H\setminus \{u,v\}$.
Let $K', N'\ge 1$ be integers, and let $0<\sigma, \sigma',\lambda',c\le 1$ with $\sigma<\sigma'<c$, such that 
$(N',K',\sigma',\lambda',c)$ forces $H'$. Then there exist integers $N,K\ge 1$
and $\lambda>0$ such that
$(N,K,\sigma,\lambda,c)$
forces $H$ $u$-first and $v$-last.
\end{thm}
\Proof Choose $\sigma''$ with $\sigma<\sigma''<\sigma'$. Choose $\lambda''>0$ and integers $K'',N''>0$ such that 
\ref{fullleafcover} is satisfied with  $\sigma, \lambda, K, K',N$ replaced by
$\sigma'', \lambda'',K'', K'+1,N''$ respectively.
Choose $\lambda>0$ and $K, N>0$ such that \ref{fullleafcover} is satisfied with
$\sigma',\lambda', K'$ replaced by $\sigma'', \lambda'', K''+1$ respectively. We may choose $N\ge N''$. 
We claim that $N,K,\sigma,\lambda$ satisfy the theorem.

Let $G$ be a $(|G|^{1-c},|G|^{1-c})$-coherent graph with $|G|\ge N$,
and let $\mathcal{A}=(A_i:i\in I)$ be a blockade in $G$ of length $K$, with shrinkage at most $\sigma$ and linkage at most $\lambda$.
We must show that there is a $u$-first and $v$-last $\mathcal{A}$-rainbow copy of $H$.
Since \ref{fullleafcover} is satisfied with
$\sigma',\lambda', K'$ replaced by $\sigma'', \lambda'', K''+1$ respectively, it follows that there exists $I''\subseteq I$ 
with $|I''|=K''+1$, such that for every partition $(P,Q)$ of $I''$,
there exists $B_h\subseteq A_h$ for each $h\in P$,
where
\begin{itemize}
\item $(B_h:h\in P)$ has shrinkage at most $\sigma''$ and linkage at most $\lambda''$; and
\item for all $h\in P$ and
all $j\in Q$ there exists $X\subseteq A_j$ that covers $B_h$ and is anticomplete to $B_{i}$ for all $i\in P\setminus \{h\}$.
\end{itemize}
Let $i_1$ be the smallest member of $I''$, and let $I_1:=\{i_1\}$, and $I_2:=I''\setminus I_1$. Choose $B_h\subseteq A_h$ 
for each $h\in I_2$ as above.

The blockade $(B_h:h\in I_2)$ has length $K''$, and shrinkage at most $\sigma''$ and linkage at most $\lambda''$. Since 
\ref{fullleafcover} is satisfied with $\sigma, \lambda, K, K',N$ replaced by
$\sigma'', \lambda'',K'', K'+1,N''$ respectively, it follows that there exists $I_3\subseteq I_2$ with $|I_3|=K'+1$, such that for every 
partition $(I_4,I_5)$ of $I_3$,
there exists $C_h\subseteq B_h$ for each $h\in I_4$,
where
\begin{itemize}
\item $(C_h:h\in I_4)$ has shrinkage at most $\sigma'$ and linkage at most $\lambda'$; and
\item for all $h\in I_4$ and
all $j\in I_5$ there exists $X\subseteq B_j$ that covers $C_h$ and is anticomplete to $C_{i}$ for all $i\in I_4\setminus \{h\}$.
\end{itemize}
Let $i_5$ be the largest member of $I_3$, and let $I_5:=\{i_5\}$ and $I_4:=I_3\setminus \{i_5\}$; and choose $(C_h:h\in I_4)$ as above.

Since $\mathcal{C}=(C_h:h\in I_4)$ has length $K'$, and shrinkage at most $\sigma'$ and linkage at most $\lambda'$, 
and $(N',K',\sigma',\lambda',c)$ forces $H'$ in $G$,
there is a $\mathcal{C}$-rainbow copy $J$ of $H'$.

Let $\phi$ be an isomorphism from $H'$ to $J$. Let $u', v'$ be the neighbours in $H$ of $u,v$ respectively. Thus 
$\phi(u'), \phi(v')\in V(J)$.
Let $h\in I_4$ such that $\phi(u')\in C_h$, and let $k\in I_4$ such that $\phi(v')\in C_k$. There exists 
$Y\subseteq B_{i_5}$ that covers $C_k$ and is anticomplete to $C_{i}$ for all $i\in I_4\setminus \{k\}$, and in particular
$Y$ is anticomplete to $V(J)\setminus \{\phi(v')\}$. Choose $y\in Y$ adjacent to $\phi(v')$.
There exists $X\subseteq A_{i_1}$ that covers $B_h$ and is anticomplete to $B_{i}$ for all $i\in I_2\setminus \{h\}$, and
in particular $X$ is anticomplete to $(V(J)\setminus \{\phi(v')\})\cup \{y\}$. Choose $x\in X$ adjacent to $\phi(u')$.
Then the subgraph induced on $V(J)\cup \{x,y\}$ is a $u$-first, $v$-last $\mathcal{A}$-rainbow copy of $H$. This proves \ref{twoleaves}.~\bbox

By combining this with \ref{addhandle}, we obtain:
\begin{thm}\label{addrealhandle}
Let $H$ be a graph, obtained from a graph $H'$ by adding a handle $P$ of length $\ell+2$.
Let $K',N'\ge 1$ be integers, and let $1\ge c,\sigma,\sigma',\lambda>0$, such that $\sigma<\sigma'$ and $c-\sigma> 1/\lfloor(\ell-1)/3\rfloor$,
and $(N',K',\sigma',\lambda',c)$ forces $H'$.
Then there exist integers $N,K>0$ and $\lambda>0$ such that
$(N,K,\sigma,\lambda,c)$ forces $H$.
\end{thm}
\Proof Let $P$ have ends $u',v'$, and let
$u,v$ be the neighbours of $u', v'$ respectively in $P$. Let $H''$ be obtained from $H$ by deleting all vertices of $P$
except $u,v,u', v'$. Thus $u,v$ are both leaves of $H''$. Choose $\sigma''$ with $\sigma'>\sigma''>\sigma$.
Choose integers $N'', K''\ge 1$ and $\lambda''>0$ such that 
\ref{twoleaves} is satisfied with $N,K',\sigma', \lambda'$ replaced by $N'',K'', \sigma'', \lambda''$ respectively.
Now $H$ is obtained from $H''$ by adding a handle of length $\ell$ with ends $u,v$. Choose $N,K,\lambda$ such that
\ref{addhandle} is satisfied with $K',\sigma', \lambda'$ replaced by $K'', \sigma'', \lambda''$ respectively. We may assume that
$N\ge N''$. We claim that $N,K,\lambda$ satisfy the theorem.

Since \ref{twoleaves} is satisfied with $N,K',\sigma', \lambda'$ replaced by $N'',K'', \sigma'', \lambda''$ respectively, and
$(N',K',\sigma',\lambda',c)$ forces $H'$, it follows that $(N,K'',\sigma'',\lambda'',c)$
forces $H$ $u$-first and $v$-last. Since \ref{addhandle} is satisfied with $K',\sigma', \lambda'$ replaced by 
$K'', \sigma'', \lambda''$ respectively, it follows that $(N,K,\sigma,\lambda,c)$ forces $H$. This proves \ref{addrealhandle}.~\bbox

\section{$\beta$-buildable graphs}

By applying \ref{addrealhandle} to $\beta$-buildable graphs, we obtain:

\begin{thm}\label{buildableforce}
Let $\beta\ge 2$ be an integer, and let $1\ge c>1/\lfloor(\beta-3)/3\rfloor$.
For every $\beta$-buildable graph $H$, and all $\sigma>0$ with $c-\sigma> 1/\lfloor(\beta-3)/3\rfloor$.
there exist integers $K,N>0$, and $\lambda>0$ such that
$(N,K,\sigma,\lambda,c)$ forces $H$.
\end{thm}
\Proof
We proceed by induction on $|H|$. If $|H|\le 2$ the result is true, so we may assume that $H$ is obtained from a smaller
$\beta$-buildable graph $H'$ by adding a handle of length at least $\beta$. Let $\ell:=\beta-2$; then 
$c-\sigma> 1/\lfloor(\ell-1)/3\rfloor$.  Choose $\sigma'$ with $\sigma<\sigma'< c-1/\lfloor(\ell-1)/3\rfloor$.
From the inductive hypothesis there exist integers $K',N'>0$, and $\lambda'>0$ such that
$(N',K',\sigma',\lambda',c)$ forces $H'$. By \ref{addrealhandle}, there exist integers $N,K>0$ and $\lambda>0$ such that
$(N,K,\sigma,\lambda,c)$ forces $H$. This proves \ref{buildableforce}.~\bbox

Now we can prove \ref{handlethm}, which we restate:
\begin{thm}\label{handlethm2}
Let $\beta\ge 2$ be an integer, let $H$ be a $\beta$-buildable graph, and let $1\ge c>1/\lfloor(\beta-3)/3\rfloor$.
There
exists $\vare>0$ such that
every $\vare$-sparse $(\vare|G|^{1-c}, \vare|G|^{1-c})$-coherent graph $G$ with $|G|>1$ contains $H$.
\end{thm}
\Proof Choose $\sigma>0$ with $c-\sigma> 1/\lfloor(\beta-3)/3\rfloor$; then by \ref{buildableforce} we can choose $K,N,\lambda$ 
such that $(N,K,\sigma,\lambda,c)$ forces $H$.
Choose $\vare>0$ such that $\vare\le 1/K$, $\vare\le 1/N$, $\vare^{\sigma}\le 1/(2K)$, and $\vare\le \lambda/(2K)$.
We claim that $\vare$ satisfies \ref{handlethm2}. 

Let $G$ be an $\vare$-sparse $(\vare|G|^{1-c}, \vare|G|^{1-c})$-coherent graph with $|G|>1$.
By \ref{big}, $|G|>1/\vare\ge K$, and so $\lfloor |G|/K\rfloor \ge |G|/(2K)$.
Consequently 
there is a blockade $\mathcal{B}$ in $G$ of length $K$ and width at least $|G|/(2K)$, and therefore with shrinkage at most $\sigma$,
since 
$|G|>1/\vare$ and $1/(2K)\ge \vare^{\sigma}\ge |G|^{-\sigma}$.

If $\mathcal{B}$ has linkage at least $\lambda$, there exist distinct $h,j\in I$, and $v\in A_h$, such that $v$ has at least $\lambda|A_j|$ neighbours in~$A_j$.
But $\lambda |A_j|\ge \lambda|G|/(2K)\ge \vare|G|$, contradicting that $G$ is $\vare$-sparse.
Hence $\mathcal{B}$ has length $K$, shrinkage at most $\sigma$ and linkage at most $\lambda$; and since $|G|>1/\vare\ge N$, and
$(N,K,\sigma,\lambda,c)$ forces $H$, there is a $\mathcal{B}$-rainbow copy of $H$. Consequently $G$ contains $H$. This proves \ref{handlethm2}.~\bbox

\section*{Acknowledgement}
Our thanks
to both referees, and especially to one who wrote an extraordinarily
thorough and helpful report.


\begin{thebibliography}{99}
\bibitem{pure1} M. Chudnovsky, A. Scott, P. Seymour and S. Spirkl, ``Pure pairs. I.  Trees and linear anticomplete pairs'',
{\em Advances in Math.}, {\bf 375} (2020), 107396, {\tt arXiv:1809.00919}, https://doi.org/10.1016/j.aim.2020.107396.
\bibitem{EHP} P. Erd\H{o}s, A. Hajnal and J. Pach, ``A Ramsey-type theorem for bipartite graphs'', {\em Geombinatorics}
  10 (2000), 64--68.
\bibitem{rodl}  V. R\"odl, ``On universality of graphs with uniformly distributed edges'',
{\em Discrete Math.} {\bf 59} (1986), 125--134.
\bibitem{pure5} A. Scott, P. Seymour and S. Spirkl, ``Pure pairs. V. Excluding some long subdivision'', 
{\em Combinatorica} {\bf 43} (2023), 571--593, {\tt arXiv:2105.03956},
https://doi.org/10.1007/s00493-023-00025-8.

\bibitem{pure6} A. Scott, P. Seymour and S. Spirkl, ``Pure pairs. VI. Excluding an
ordered tree'', {\em SIAM J. Disc. Math.}, {\bf 36} (2022), 170--187, {\tt
arXiv:2009.10671}, https://doi.org/10.1137/20M1368331.

\end{thebibliography}
\end{document}